\documentclass[12pt]{amsart}

\usepackage{amssymb}
\usepackage{amsthm}
\usepackage{amscd}

\usepackage{amsmath}
\usepackage{epic,eepic}
\usepackage{graphicx}

\usepackage{latexsym}

\usepackage{pifont}

\theoremstyle{plain}
\newtheorem{lemma}{Lemma}[subsection]
\newtheorem{prop}[lemma]{Proposition}
\newtheorem{thm}[lemma]{Theorem}
\newtheorem{cor}[lemma]{Corollary}
\newtheorem{aplemma}{Lemma~A.\hspace{-1.5mm}}
\newtheorem{approp}{Proposition~A.\hspace{-1.5mm}}
\newtheorem{apthm}{Theorem~A.\hspace{-1.5mm}}
\newtheorem{apcor}{Corollary~A.\hspace{-1.5mm}}
\newtheorem{intthm}{Theorem}

\usepackage[all]{xy}

\newtheorem{conj}[lemma]{Conjecture}

\theoremstyle{definition}

\newtheorem{rema}[lemma]{Remark}

\newtheorem{remb}{Remark}

\newtheorem{defi}[lemma]{Definition}
\newtheorem{exa}[lemma]{Example}
\newtheorem{aprem}{Remark~A.\hspace{-1.5mm}}
\newtheorem{apdefi}{Definition~A.\hspace{-1.5mm}}
\newcommand{\bde}{\begin{defi}}
\newcommand{\ede}{\end{defi}\vspace{1mm}}
\newcommand{\ble}{\begin{lemma}}
\newcommand{\ele}{\end{lemma}}
\newcommand{\bpr}{\begin{prop}}
\newcommand{\epr}{\end{prop}}
\newcommand{\bt}{\begin{thm}}
\newcommand{\et}{\end{thm}}
\newcommand{\bco}{\begin{cor}}
\newcommand{\eco}{\end{cor}}
\newcommand{\bre}{\begin{rema}}
\newcommand{\ere}{\end{rema}}
\newcommand{\brea}{\begin{rema}}
\newcommand{\erea}{\end{rema}\vspace{1mm}}
\newcommand{\breb}{\begin{remb}}
\newcommand{\ereb}{\end{remb}\vspace{1mm}}
\newcommand{\bex}{\begin{exa}}
\newcommand{\eex}{\end{exa}}
\newcommand{\bpf}{\begin{proof}}
\newcommand{\epf}{\end{proof}\vspace{1mm}}

\newcommand{\bade}{\begin{apdefi}}
\newcommand{\eade}{\end{apdefi}}
\newcommand{\bale}{\begin{aplemma}}
\newcommand{\eale}{\end{aplemma}}
\newcommand{\bapr}{\begin{approp}}
\newcommand{\eapr}{\end{approp}}
\newcommand{\bat}{\begin{apthm}}
\newcommand{\eat}{\end{apthm}}
\newcommand{\baco}{\begin{apcor}}
\newcommand{\eaco}{\end{apcor}}
\newcommand{\bare}{\begin{aprem}}
\newcommand{\eare}{\end{aprem}}

\newcommand{\be}{\begin{enumerate}}
\newcommand{\ee}{\end{enumerate}}
\newcommand{\bcd}{\[\begin{CD}}
\newcommand{\ecd}{\end{CD}\]}
\newcommand{\bit}{\begin{itemize}}
\newcommand{\eit}{\end{itemize}}
\newcommand{\bq}{\begin{quote}}
\newcommand{\eq}{\end{quote}}
\newcommand{\ba}{\begin{array}}
\newcommand{\ea}{\end{array}}

\newcommand{\mcD}{\mathcal{D}}
\newcommand{\mcE}{\mathcal{E}}

\newcommand{\mcO}{\mathcal{O}}

\newcommand{\mcT}{\mathcal{T}}

\newcommand{\mbC}{\mathbb{C}}

\newcommand{\mbF}{\mathbb{F}}

\newcommand{\mbP}{\mathbb{P}}

\newcommand{\mbZ}{\mathbb{Z}}

\newcommand{\mfB}{\mathfrak{B}}
\newcommand{\mfC}{\mathfrak{C}}

\newcommand{\mfE}{\mathfrak{E}}

\newcommand{\mfM}{\mathfrak{M}}

\newcommand{\mfO}{\mathfrak{O}}
\newcommand{\mfP}{\mathfrak{P}}

\newcommand{\mfT}{\mathfrak{T}}

\newcommand{\mfX}{\mathfrak{X}}

\newcommand{\mfa}{\mathfrak{a}}

\newcommand{\mfc}{\mathfrak{c}}
\newcommand{\mfd}{\mathfrak{d}}
\newcommand{\mfe}{\mathfrak{e}}

\newcommand{\mfl}{\mathfrak{l}}

\newcommand{\mfn}{\mathfrak{n}}
\newcommand{\mfo}{\mathfrak{o}}
\newcommand{\mfp}{\mathfrak{p}}

\newcommand{\mfr}{\mathfrak{r}}
\newcommand{\mfs}{\mathfrak{s}}
\newcommand{\mft}{\mathfrak{t}}
\newcommand{\mfu}{\mathfrak{u}}

\newcommand{\migi}{\rightarrow}

\newcommand{\isom}{\stackrel{\sim}{\migi}}

\newcommand{\migiincl}{\hookrightarrow}

\newcommand{\migisurj}{\twoheadrightarrow}

\newcommand{\mr}{\mathrm}
\newcommand{\hidden}[1]{\,}

\DeclareMathAlphabet{\mathpzc}{OT1}{pzc}{m}{it}
\newcommand{\vin}{\rotatebox{90}{$\in$}}

\begin{document}

\title[Dormant Miura Transformation]{A Combinatorial Description of  \\ The Dormant
Miura Transformation}
\author{Yasuhiro Wakabayashi}
\date{}
\markboth{Yasuhiro Wakabayashi}{}
\maketitle
\footnotetext{Y. Wakabayashi: Department of Mathematics, Tokyo Institute of Technology, 2-12-1 Ookayama, Meguro-ku, Tokyo 152-8551, JAPAN;}
\footnotetext{e-mail: {\tt wkbysh@math.titech.ac.jp};}
\begin{abstract}
A dormant generic Miura $\mathfrak{sl}_2$-oper is a flat $\mathrm{PGL}_2$-bundle over an algebraic curve in positive characteristic equipped with some additional data. In the present paper, we give a combinatorial description of dormant generic Miura $\mfs \mfl_2$-opers on a totally degenerate curve. The combinatorial objects that we use are certain branch numberings of $3$-regular graphs. Our description may be thought of as an analogue of the combinatorial description of dormant $\mfs \mfl_2$-opers given by S. Mochizuki, F. Liu, and B. Osserman. It allows us to think of the Miura transformation in terms of combinatorics. As an application, we identify the dormant generic Miura $\mfs \mfl_2$-opers on totally degenerate curves of genus $>0$.
\end{abstract}
\tableofcontents 

\section*{Introduction}

\vspace{5mm}
\subsection*{0.1} \label{S01}

The purpose of the present paper  is to give a combinatorial description of  
dormant generic Miura $\mfs \mfl_2$-opers.
The combinatorial objects that we use are certain branch numberings of $3$-regular (i.e., trivalent) graphs.
Our description will  be thought of as an analogue of the combinatorial description of dormant $\mfs \mfl_2$-opers given by S. Mochizuki, F. Liu, and B. Osserman (cf. ~\cite{Mzk2}, Introduction, \S\,1.2, Theorem 1.3; ~\cite{LO}, Theorem 3.9).
It allows us to think of the {\it Miura transformation} in terms of combinatorics.

The celebrated {\it Miura transformation} concerns the Korteweg-de Vries (KdV) and the modified KdV (mKdV) equations.
The KdV equation was derived as an equation modeling the behavior of shallow water waves moving in one direction by Korteweg and his student de Vries.
The {\it KdV equation} reads 
\begin{align}
v_t - 6 v v_x + v_{xxx} =0,
\end{align}
while its modified counterpart, the {\it mKdV equation}, equals
\begin{align}
u_t - 6 u^2 u_x + u_{xxx} =0.
\end{align}
Let  us take a solution $u \in \mbC [[t, x]]$ to the mKdV equation.
Then, the function $v$  characterized  by  the equality of differential operators
\begin{align} \label{EE010}
(\partial_x - u) (\partial_x + u) = \partial^2_x - v  
\end{align}
(i.e., $v  = u^2 - u_x$) specifies a solution to  the KdV equation.
The assignment  
\begin{align} \label{EE012}
u \mapsto  v  \ (= u^2 - u_x)
\end{align}
 is nothing but  the  Miura transformation.

Recall that the differential operator $\partial_x^2 - v$ (resp., $(\partial_x - u)(\partial_x + u)$) in the right-hand (resp., the left-hand) side of (\ref{EE010}) corresponds, in the usual manner,  to 
the connection on a vector bundle (i.e.,
 the first order matrix differential operator)  of the form
\begin{align} \label{EE014}
\nabla = \partial_x - \begin{pmatrix} 0 & v \\ 1& 0   \end{pmatrix} \ \left(\text{resp.,} \ \nabla_{\mr{Miura}} = \partial_x - \begin{pmatrix} u & 0 \\ 1& -u   \end{pmatrix}\right).
\end{align}
A connection  of the form (\ref{EE014}) is called 
an $\mfs \mfl_2$-oper (resp., a generic Miura $\mfs \mfl_2$-oper).
If  the underlying space of the vector bundle  is a Riemann surface, then such a  connection may be identified with  a projective connection   (resp., an affine connection) on this Riemann surface (cf. ~\cite{G2}).
$\nabla_{\mr{Miura}}$ becomes  $\nabla$ after  gauge transformation by some upper triangular matrix.
This implies that the Miura transformation may be thought of as 
the assignment from generic Miura $\mfs \mfl_2$-opers $\nabla_{\mr{Miura}}$ to $\mfs \mfl_2$-opers $\nabla$ induced  by gauge transformations  in this way.

\vspace{5mm}
\subsection*{0.2} \label{S02}

In the present paper, we consider the case of (Miura) opers {\it in positive characteristic}.
A remarkable point is that
unlike the complex case,
 one may have a generic Miura $\mfs \mfl_2$-oper on some  entire (i.e., proper) smooth curve in positive characteristic of genus $>1$.
This fact has  already proved in the previous work concerning {\it dormant} generic Miura opers and Tango structures  (cf. ~\cite{Wak7},  Theorems A and B).
Here, an $\mfs \mfl_2$-oper (resp., a generic Miura $\mfs \mfl_2$-oper) is called  dormant if it has vanishing $p$-curvature.
(We refer to ~\cite{Wak5} for various discussions and results concerning  dormant opers on pointed stable curves.)
Each pointed stable curve $\mfX$ gives rise to 
the set
\begin{align}
\overline{\mfO} \mfp_{\mfs \mfl_2, \mfX}^{^\mr{Zzz...}} \ \  \ (\text{resp.,} \  \mfM \overline{\mfO} \mfp_{\mfs \mfl_2, \mfX}^{^\mr{Zzz...}})
\end{align}
of isomorphism classes of dormant $\mfs \mfl_2$-opers (resp., dormant generic Miura $\mfs \mfl_2$-opers) on $\mfX$.
The Miura transformation described in terms of opers (i.e., $\nabla_{\mr{Miura}} \mapsto \nabla$ as discussed in the previous subsection) gives  
 a map
\begin{align} \label{AAAA}
\mu^{^\mr{Zzz...}}_\mfX : \mfM \overline{\mfO} \mfp_{\mfs \mfl_2,  \mfX}^{^\mr{Zzz...}} \migi \overline{\mfO} \mfp_{\mfs \mfl_2, \mfX}^{^\mr{Zzz...}},
\end{align}
referred  to as the {\bf dormant Miura transformation}.
Our main interest of the present paper is to understand,  in terms of combinatorics,  the behavior of 
this map for the case where the underlying curve is totally degenerate (cf. the discussion preceding Proposition  \ref{y0176})

 \vspace{5mm}
\subsection*{0.3} \label{S03}
 
To this end, we first recall the  combinatorial description of dormant $\mfs \mfl_2$-opers on such a curve, which was obtained  in the work 
by S. Mochizuki  in the context of
   $p$-adic Teichm\"{u}ller theory
    and the work by F. Liu and B. Osserman.
Let  $(g, r)$ be a pair of nonnegative integers with $2g-2+r >0$ and 
$\mfX$
  a  totally degenerate curve  of type $(g,r)$  over an algebraically closed field $k$ of characteristic $p >2$.
$\mfX$ associates its dual marked semi-graph
$\Gamma^+$
  (cf. the discussion in Remark \ref{R00f1}), which is  $3$-regular.
A {\bf balanced $p$-edge numbering} (cf. Definition \ref{De113}, (i) for its precise definition) on $\Gamma^+$ is 
a numbering $(a_e)_{e \in E_\Gamma}$ on the set of  edges $E_\Gamma$ 
of $\Gamma^+$
 satisfying certain triangle inequalities with respect to each triple of numbers incident to one vertex.
Denote 
by 
\begin{align}
p \text{-}\mfE \mfd_{\Gamma^+}^\mr{bal}
\end{align}
 the set of balanced $p$-branch numberings  on $\Gamma^+$. 
Then, there exists a canonical  bijection
\begin{align} \label{EE016}
p\text{-}\mfE \mfd_{\Gamma^+}^{\mr{bal}} \isom \overline{\mfO} \mfp_{\mfs \mfl_2, \mfX}^{^\mr{Zzz...}} 
\end{align}
(cf. the discussion preceding Proposition \ref{y0176} for the precise construction).
The  inverse of (\ref{EE016}) is obtained by taking the {\it radii} of  dormant $\mfs \mfl_2$-opers on $\mfX$ restricted to the various  rational  components of its normalization.

\vspace{5mm}

\hspace{10mm}
\begin{picture}(200,100)
\put(30,-15){{\bf A balanced $p$-edge}}
\put(25,-25){{\bf numbering (with $p\,\ooalign{$>$ \cr $\ \, >$}\hspace{1mm}0$)}}

\qbezier(10,40)(35,110)(90,70)
\qbezier(10,40)(50,60)(90,70)
\qbezier(10,40)(40,20)(70,30)
\qbezier(70,30)(80,40)(90,70)
\qbezier(70,30)(120,0)(135,45)
\put(140,64){\circle{40}}
\put(10,40){\circle*{5}}
\put(90,70){\circle*{5}}
\put(70,30){\circle*{5}}
\put(135,45){\circle*{5}}
\put(123,25){\circle*{5}}
\put(163,15){\circle*{5}}
\put(140,7){$3$}
\put(132,30){$4$}
\qbezier(163,15)(143,15)(123,25)
\put(13,65){$6$}
\put(40,60){$7$}
\put(32,17){$5$}
\put(81,40){$3$}
\put(85,7){$2$}
\put(145,50){$2$}


\end{picture}
\hspace{10mm}
\begin{picture}(200,100)
\put(5,40){\circle*{5}}
\put(47,37){\circle*{5}}
\put(60,5){\circle*{5}}
\put(71,70){\circle*{5}}
\put(45,90){\circle*{5}}
\put(110,68){\circle*{5}}
\put(130,90){\circle*{5}}
\put(130,30){\circle*{5}}

\qbezier(5,40)(26,38)(47,37)
\qbezier(60,5)(53,21)(47,37)
\qbezier(71,70)(59,53)(47,37)
\qbezier(71,70)(58,80)(45,90)
\qbezier(71,70)(90,69)(110,68)
\qbezier(130,90)(120,79)(110,68)
\qbezier(130,30)(120,49)(110,68)

\put(62,8){$9$}
\put(5,30){$10$}
\put(33,28){$1$}
\put(52,25){$2$}
\put(47,48){$9$}
\put(67,56){$2$}
\put(57,69){$3$}
\put(50,88){$8$}
\put(76,73){$7$}
\put(100,71){$4$}
\put(117,68){$5$}
\put(120,87){$6$}
\put(119,55){$3$}
\put(135,30){$8$}
\put(30,-15){{\bf A strict  $p$-branch}}
\put(25,-25){{\bf numbering (with $p =11$)}}
\end{picture}

\vspace{15mm}

 Next, we introduce the notion of  a   {\bf  strict $p$-branch numbering} (cf. Definition \ref{FDE2}) on $\Gamma^+$, which  is defined to be 
  a certain numbering $(\varepsilon_b)_{b \in B_{\Gamma}}$ on the set of branches $B_\Gamma$ of $\Gamma^+$ satisfying some conditions, including  the  condition  that 
  the three numbers incident  to one   vertex amount precisely to $p+1$.
Each strict $p$-branch numbering  $(\varepsilon_b)_{b \in B_{\Gamma}}$ on $\Gamma^+$
induces  a balanced $p$-edge numbering $(\varepsilon^{\mu}_e)_{e \in E_{\Gamma}}$ on $\Gamma^+$ which is well-defined in such a way that
for each  $e \in E_{\Gamma}$ having  a branch $b \in B_{\Gamma}$, 
$\varepsilon^{\mu}_e$ equals $\frac{p-\varepsilon_e -1}{2}$ if $\varepsilon_e$ is even and $\varepsilon^{\mu}_e$ equals $\frac{m-1}{2}$ if $\varepsilon_e$ is odd.
by
\begin{align}  
  p\text{-}\mfB  \mfr_{\Gamma^+}^\mr{st}
\end{align}
  the set of  strict    $p$-branch numberings on $\Gamma^+$.
  Then,  
the assignment $(\varepsilon_b)_{b \in B_{\Gamma}} \mapsto (\varepsilon^{\mu}_{e})_{e \in E_{\Gamma}}$ defines 
 a map 
 \begin{align} \label{eq4730}
\mu^\mr{comb}_{\Gamma^+} : p\text{-}\mfB  \mfr_{\Gamma^+}^\mr{st} & \migi p\text{-}\mfE \mfd^\mr{bal}_{\Gamma^+}. 
 \end{align}
The main result of the present paper is as follows (cf. Corollary \ref{PP005} for its refinement).

\vspace{3mm}
\begin{intthm} \label{ThB} 
\leavevmode\\
 \ \ \ 
 Let $\mfX$ and $\Gamma^+$ be as above.
 Then, there exists a canonical bijection
 \begin{align} \label{EErh9}
p\text{-} \mfB \mfr_{\Gamma^+}^\mr{st} \isom  \mfM \overline{\mfO} \mfp_{\mfs \mfl_2, \mfX}^{^\mr{Zzz...}}   
\end{align}
 making  the following square diagram commute:
\begin{align} \label{RR2}
\begin{CD}
 p\text{-} \mfB \mfr_{\Gamma^+}^\mr{st} @> (\ref{EErh9}) > \sim > \mfM \overline{\mfO} \mfp_{\mfs \mfl_2, \mfX}^{^\mr{Zzz...}} 
\\
@V \mu^\mr{comb}_{\Gamma^+}  VV @VV \mu^{^\mr{Zzz...}}_{\mfX}  V
\\
 p\text{-} \mfE \mfd_{\Gamma^+}^\mr{bal} @> \sim > (\ref{EE016}) > \overline{\mfO} \mfp^{^\mr{Zzz...}}_{\mfs \mfl_2, \mfX}.
\end{CD}
\end{align}

 \end{intthm}
\vspace{3mm}

The above theorem implies that  one may consider the morphism $\mu^{\mr{comb}}_{\Gamma^+}$ as a  combinatorial realization of the dormant Miura transformation.
By means of this combinatorial description,
 we  prove (cf. Corollary \ref{PCD04}) that
there is no dormant generic Miura $\mfs \mfl_2$-oper on a  totally degenerate curve of genus $>1$. 
This result  may be  thought of as a combinatorial (and dormant) analogue of the emptiness (proved in ~\cite{G2},  Lemma 1) of the space of complex affine structures on a compact hyperbolic Riemann surface.
Also, we identify (in terms of combinatorics) the dormant generic Miura $\mfs \mfl_2$-opers on totally degenerate curves of genus $1$ (cf. Proposition \ref{P048} (ii) and its proof).

\vspace{5mm}
\hspace{-4mm}{\bf Acknowledgement} \leavevmode\\
 \ \ \ 
The author cannot express enough his sincere and  deep gratitude to  all those who give  the opportunity or   impart  the great joy of  studying 
  mathematics to him.
The author wrote  the present paper for a gratitude letter to them.
The author was partially  supported by the Grant-in-Aid for Scientific Research (KAKENHI No.\,18K13385).


\vspace{5mm}
\hspace{-4mm}{\bf Notation and conventions} \leavevmode\\
 \ \ \ 
Let us introduce some  notation and conventions used in  the present paper.
Throughout the present paper, we  fix an algebraically closed field  $k$ of characteristic $p>2$ and  a pair of nonnegative integers $(g,r)$ with $2g-2 +r >0$. 

For a log  scheme indicated, say,  by $Y^\mr{log}$, we shall write $Y$ for the underlying scheme of $Y^\mr{log}$.
If, moreover,  
 $Z^\mr{log}$ is a  log scheme  over $Y^\mr{log}$,
then we shall write
$\Omega_{Z^\mr{log}/Y^\mr{log}}$ for the sheaf of  logarithmic $1$-forms on $Z^\mr{log}$ over $Y^\mr{log}$ and write $\mcT_{Z^\mr{log}/Y^\mr{log}} := \Omega_{Z^\mr{log}/Y^\mr{log}}^\vee$ for its dual. 
(Basic references for the notion of a log scheme  are ~\cite{KaKa}, ~\cite{ILL2}, and ~\cite{KaFu}.)

Given a set $S$ and a positive integer $r$, we shall denote by
$S^{\times r}$ the set of $r$-tuples of elements in $S$.




\vspace{10mm}
\section{Dormant generic Miura $\mfs \mfl_2$-opers} \vspace{3mm}


\vspace{5mm}
\subsection{Semi-graphs}
\leavevmode\\ \vspace{-4mm}

First, recall from ~\cite{Mzk3}, Appendix (or ~\cite{Mzk4}, \S\,1), the definition of a semi-graph, as follows.

\vspace{3mm}
\bde \leavevmode\\
 \  \ \ 
 A  {\bf  semi-graph}  is a triple 
\begin{align}
\Gamma = (V_\Gamma, E_\Gamma, \zeta_\Gamma),
\end{align}
 where
 \begin{itemize}
 \item[$\bullet$]
  $V_\Gamma$ denotes a set, whose elements are
 called {\bf vertices};
 \item[$\bullet$]
 $E_\Gamma$ denotes a  set, whose elements are
  called {\bf edges},  consisting of   sets  with cardinality $2$ satisfying the condition that $e  \neq e' \in E_\Gamma$ implies $e \cap e' = \emptyset$;
 \item[$\bullet$]
 $\zeta_\Gamma$ denotes a  map $\coprod_{e \in E_\Gamma} e \migi V_\Gamma \cup \{ V_\Gamma \}$ (where we note that $V_\Gamma \cap \{ V_\Gamma \} = \emptyset$ since $V_\Gamma \notin V_\Gamma$), which is called a {\bf coincidence map}.
 \end{itemize}

Let $e$ be
an edge  of $\Gamma$ (i.e., $e \in E_\Gamma$).
Then,  we shall refer to any element $b \in e$   as a {\bf branch}  of  $e$.
If $b$ is a branch of $e$, then we shall denote by 
\begin{align} \label{EE001}
b^\divideontimes \ \left(\in e \right)
\end{align}
the branch of $e$ with $\{ b, b^\divideontimes \} = e$.
We shall  write
\begin{align}
B_\Gamma := \coprod_{e \in E_\Gamma} e
\end{align} 
and, for each 
 $v\in V_\Gamma \cup \{ V_\Gamma\}$, we write
 \begin{align}
 B_v  :=  \zeta_\Gamma^{-1} (\{ v\})
 \end{align}
(hence, $B_\Gamma = \coprod_{v \in V_\Gamma \cup \{ V_\Gamma\}} B_v$).
  \vspace{3mm}
  \ede

Let us fix a semi-graph  $\Gamma := (V_\Gamma, E_\Gamma, \zeta_\Gamma)$.

\bde \leavevmode\\
\vspace{-5mm}
\begin{itemize}
\item[(i)]
We shall say that $\Gamma$ is {\bf finite} if both $V_\Gamma$ and $E_\Gamma$ are finite.
\item[(ii)]
 Let $m$ be a positive  integer.
We shall say that  $\Gamma$ is {\bf $m$-regular} if for  any vertex $v \in V_\Gamma$, the cardinality of  $B_v$ is precisely $m$.
\end{itemize}
  \ede

\bde \leavevmode\\
\vspace{-5mm}
\begin{itemize}
\item[(i)]
Let $u$, $v$ be vertices of $\Gamma$.
A {\bf path} from $u$ to $v$ is an ordered  collection $(b_j)_{j=1}^l$ (for some positive integer $l$) of branches of $\Gamma$
such that  
$u = \zeta_\Gamma (b_1)$,  $v = \zeta_\Gamma (b_l^\divideontimes)$, and  $\zeta_\Gamma (b_j^\divideontimes) = \zeta_\Gamma (b_{j+1})$  for any $j \in \{ 1, \cdots, l-1 \}$.
\item[(ii)]
We shall say that  $\Gamma$ is {\bf connected} if for any two distinct  vertices $u$, $v \in V_\Gamma$,
 there exists  
 a path from $u$ to $v$.
\end{itemize}
  \ede

\bde \leavevmode\\
 \ \ \ 
 Suppose that  $\Gamma$ is
 finite, connected, and $3$-regular. 
Then, we shall say that $\Gamma$ is {\bf of type $(g,r)$}
if  the following equalities hold:
\begin{align}
g = 1  - \sharp (V_\Gamma) + \sharp (E_\Gamma) - \sharp (B_{V_\Gamma}), \hspace{10mm}
r = \sharp (B_{V_\Gamma}).
\end{align}

  \ede

\vspace{3mm}
\begin{rema} \label{R0ff01}
\leavevmode\\
 \ \ \   In a natural way (cf., e.g.,  ~\cite{Mzk4}, \S\,1), the  semi-graph  $\Gamma$ may be thought of as a {\it topological space}.
Thus, it makes sense to speak of the {\it Betti number} of $\Gamma$, which we shall denote  by $\beta (\Gamma)$.
If $\Gamma$ is finite, connected, and $3$-regular, then  the equality $g = \beta (\Gamma)$ holds.
\end{rema}
\vspace{3mm}

\bde \leavevmode\\
\vspace{-5mm}
\begin{itemize}
\item[(i)]
Suppose that $\Gamma$ is finite.
A {\bf marking} on $\Gamma$ is a bijection of sets  
\begin{align}
\lambda_{\Gamma} : B_{V_\Gamma} \isom \{ 1, \cdots, r \},
\end{align}
 where $r := \sharp (B_{V_\Gamma})$.
In particular,  if
$\sharp (B_{\Gamma_V}) =0$, then we consider any semi-graph  as being equipped with a unique marking  $B_{\Gamma_V}\isom \emptyset$.
\item[(ii)]
A {\bf marked semi-graph} is a quadruple
\begin{align}
\Gamma^+ := (V_\Gamma, E_\Gamma, \zeta_\Gamma, \lambda_\Gamma),
\end{align}
where $\Gamma := (V_\Gamma, E_\Gamma, \zeta_\Gamma)$ is a finite semi-graph and $\lambda_\Gamma$ is a marking on $\Gamma$.
We shall refer to $\Gamma$ as the {\bf underlying semi-graph} of $\Gamma^+$.
\end{itemize} 
 \ede

\vspace{3mm}
\begin{rema} \label{R00f1}
\leavevmode\\
\ \ \ 
Let $\mfX := (X/k, \{ \sigma_i \}_{i=1}^r)$ be a pointed stable curve   of type $(g,r)$  over $k$.
Then,  in the usual manner, one can associate to $\mfX$   a  marked  semi-graph 
\begin{align} \label{EEE445}
\Gamma_{\mfX}^+ := (V_{\mfX}, E_{\mfX}, \zeta_{\mfX}, \lambda_\mfX)
\end{align}
defined as follows.
$V_{\mfX}$ is  the set of irreducible components of $X$ and  $E_{\mfX}$ is  the disjoint union $N_{\mfX} \sqcup  \{ \sigma_i \}_{i=1}^r$
 of  the set of nodal points  $N_{\mfX}$ ($\subseteq X (k)$)  and the set of marked points $\{ \sigma_i \}_{i=1}^r$.
 Here,  note that any  node $e \in X(k)$ has two distinct  branches $b_1$ and $b_2$, each of which lies on some well-defined irreducible component of $X$;
 we shall identify  $e$ 
  with   $\{ b_1, b_2\}$.  Moreover,  identify each marked point $\sigma_i$ with the set $\{ \sigma_i, \{ \sigma_i \} \}$.
 In this way, we regard the  elements of $E_{\mfX}$ as  sets  with cardinality $2$.
Also, we define $\zeta_{\mfX}$ to be  the  map  $\coprod_{e \in E_{\mfX}} e  \migi V_{\mfX} \cup \{ V_{\mfX} \}$ determined as follows:
\begin{itemize}
\item[$\bullet$]
if $b \in e$ (for some $e \in N_{\mfX}$) or $b = \sigma_i$ (for some $i$), then
$\zeta_{\mfX} (b)$ is the irreducible component  (i.e., an element of $V_{\mfX}$) in which $b$ lies;
\item[$\bullet$]
 $\zeta_{\mfX} (\{ \sigma_i \}) = \{V_{\mfX}\}$ for any $i$.
\end{itemize}
(Hence, $B_{V_\mfX} = \{ \{ \sigma_i \} \}_{i=1}^r$.)
Finally, $\lambda_\mfX$ is given by assigning $\{ \sigma_i \} \mapsto i$ (for any $i \in \{1, \cdots, r \}$).
We shall refer to $\Gamma_{\mfX}^+$ as the {\bf dual marked semi-graph} associated with $\mfX$.
\end{rema}
\vspace{3mm}


We shall write
\begin{align}
\widetilde{\mbF}_p := \{ 0,1, \cdots, p-1 \} \ \left(\subseteq \mbZ \right),
\end{align}
and write $\tau$ for the natural composite bijection
\begin{align}
\tau : \widetilde{\mbF}_p \migiincl \mbZ \migisurj \mbF_p \ \left(:= \mbZ/p \mbZ \right).
\end{align}
Let us  define an involution $(-)^\veebar$ on  $\widetilde{\mbF}_p$ to be the map  given as follows:
 \begin{equation} \label{w010w}
m^\veebar := \begin{cases} p-m & \text{if $m \in \widetilde{\mbF}_p \setminus \{ 0 \}$,} \\
 \ \ \ 0 & \text{if $m = 0$.}
\end{cases}
\end{equation}
In particular, we have $\tau (m^\veebar) = - \tau (m)$.

\vspace{3mm}
\bde \label{y0175}\leavevmode\\
\vspace{-5mm}
\begin{itemize}
\item[(i)]
A {\bf $p$-branch numbering} on $\Gamma$ is a collection $\vec{m} := (m_b)_{b \in B_\Gamma} \in \widetilde{\mbF}_p^{B_\Gamma}$ of elements of $\widetilde{\mbF}_p$ indexed by the set 
$B_\Gamma$
 such that
for any edge $e := \{ b, b^\divideontimes \}$, the equality $m_{b^\divideontimes} = m_{b}^\veebar$ holds.
\item[(ii)]
Assume  further that $\Gamma$ is finite and $B_{V_\Gamma} \neq \emptyset$.
Also, assume that  we are given a marking $\lambda_\Gamma : B_{V_\Gamma} \isom \{1, \cdots, r \}$ (where  $r := \sharp (B_{V_\Gamma})$) on $\Gamma$ and  an element $\vec{\varepsilon} := (\varepsilon_i)_{i=1}^r$ of $\mbF_p^{\times r}$.
Then, we shall say that a $p$-branch numbering 
  $\vec{m} := (m_b)_{b \in B_\Gamma}  \in \widetilde{\mbF}_p^{B_\Gamma}$ 
 is {\bf of exponent  $\vec{\varepsilon}$} if 
 $\tau (m_{b}) = \varepsilon_{\lambda_\Gamma (b)}$ for any $b \in B_{V_\Gamma}$.
For convenience, (regardless of whether $B_{V_\Gamma}$ is empty or not) we shall refer to  any $p$-branch numbering  as being {\bf of exponent  $\emptyset$}.
\end{itemize}
 \ede

\vspace{5mm}
\subsection{Opers and Miura opers on a pointed stable curve} \label{z041}
\leavevmode\\ \vspace{-4mm}

Next,  we recall
the notion of a dormant generic Miura oper defined on a pointed stable curve.
We refer to ~\cite{Wak7} for various definitions and notation used in this section.

Let 
\begin{equation}
\label{X}
 \mfX : =(f :X \migi \mr{Spec} (k), \{ \sigma_i : \mr{Spec} (k) \migi X\}_{i=1}^r)\end{equation}
be  a pointed stable curve over $k$ of type $(g,r)$, consisting of a (proper) semi-stable curve $X$ over $k$ of genus $g$ and $r$ marked points $\sigma_i$ ($i = 1, \cdots, r$) of $X$.
Note that there exists  natural log structures on $X$ and $\mr{Spec} (k)$ (cf. ~\cite{Wak7}, \S\,1.1);
 we denote  the resulting log schemes by 
\begin{align}
\mr{Spec}(k)^{\mfX \text{-}\mr{log}}  \ \ \text{and}  \ \ X^{\mfX \text{-}\mr{log}}
\end{align} 
respectively.
If there is no fear of causing confusion, we write $\mr{Spec}(k)^\mr{log}$ (or just $k^\mr{log}$) and $X^\mr{log}$ instead of $\mr{Spec} (k)^{\mfX \text{-}\mr{log}}$ and  $X^{\mfX \text{-}\mr{log}}$ respectively.
The structure morphism $f : X \migi \mr{Spec}(k)$ extends to a morphism $f^\mr{log} : X^{\mr{log}} \migi \mr{Spec}(k)^{\mr{log}}$ of log schemes, 
by which $X^\mr{log}$ determines a log-curve over $\mr{Spec}(k)^\mr{log}$ 
(cf.  ~\cite{ACGH},  Definition 4.5 for the definition of a log-curve).

Let us recall briefly the definitions of an $\mfs \mfl_2$-opers and a Miura $\mfs \mfl_2$-oper.
 Denote by $\mr{PGL}_2$  the projective linear group over $\mbF_p$ of rank $2$ (considered as an algebraic group over $k$ via base-change by $\mbF_p \migiincl k$) and by $B$  the Borel subgroup of $\mr{PGL}_2$ defined to be the image (via the quotient $\mr{GL}_2 \migisurj \mr{PGL}_2$) of upper triangular matrices.
Let us identify  $\mfs \mfl_2$ with  the Lie algebra of $\mr{PGL}_2$.
 An {\bf $\mfs \mfl_2$-oper} 
  on $\mfX$ is a pair $\mcE^\spadesuit := (\mcE_B, \nabla_\mcE)$ consisting of a (right) $B$-torsor $\mcE_B$ over $X$ and a $k^\mr{log}$-connection $\nabla_\mcE$ 
on the $\mr{PGL}_2$-torsor $\mcE_{\mr{PGL}_2} := \mcE_B \times^B \mr{PGL}_2$ induced by $\mcE_B$  such that $\mcE_B$ is transversal to $\nabla_\mcE$.
 (Here, we refer to ~\cite{Wak7}, \S\,1.3,  for the definition of a connection on a torsor  in the  logarithmic sense,  and refer to   Definition 3.1.1 (i) in {\it loc.\,cit.} or ~\cite{Mzk2}, Chap.\,I, \S\,2, Definition 2.2 for the  precise definition of an $\mfs \mfl_2$-oper).
Also, a {\bf  Miura $\mfs \mfl_2$-oper} on $\mfX$ is defined to be a collection of data
$\widehat{\mcE}^\spadesuit := (\mcE_B, \nabla_\mcE, \mcE'_B, \eta_\mcE)$, where
$(\mcE_B, \nabla_\mcE)$ is  an $\mfs \mfl_2$-oper  on $\mfX$,  
 $\mcE'_B$ is another  $B$-torsor $\mcE'_B$ over $X$, and  $\eta_\mcE$ is 
an isomorphism $\mcE'_B \times^B \mr{PGL}_2 \isom \mcE_{\mr{PGL}_2}$ of $\mr{PGL}_2$-torsors via which   $\mcE'_B$ is preserved by $\nabla_\mcE$   (cf. Definition 3.2.1 in ~\cite{Wak7} for the precise definition of a Miura $\mfs \mfl_2$-oper). 
We shall say that a Miura $\mfs \mfl_2$-oper $\widehat{\mcE}^\spadesuit := (\mcE_B, \nabla_\mcE, \mcE'_B, \eta_\mcE)$ is {\bf generic}   (cf. Definition 3.3.1 in {\it loc.\,cit.}) if $\mcE_B$ and $\mcE'_B$ are in generic relative position.
Moreover, we shall say that an $\mfs \mfl_2$-oper $\mcE^\spadesuit := (\mcE_B,\nabla_\mcE)$ (resp., a  Miura $\mfs \mfl_2$-oper $\widehat{\mcE}^\spadesuit := (\mcE_B, \nabla_\mcE, \mcE'_B, \eta_\mcE)$) is {\bf dormant} if  $\nabla_\mcE$ has vanishing $p$-curvature (cf. ~\cite{Wak5}, Definition 3.2.1 for the definition of $p$-curvature).
Denote by 
\begin{align} \label{eqeq100}
 \overline{\mfO} \mfp_{\mfs \mfl_2, \mfX} \ \ \ \left(\text{resp.,} \ \mfM \overline{\mfO} \mfp_{\mfs \mfl_2, \mfX}\right)
\end{align}
the set of isomorphism classes of  $\mfs \mfl_2$-opers (resp., the set of isomorphism classes of generic Miura $\mfs \mfl_2$-opers) on $\mfX$.
Also,
denote by 
\begin{align} \label{eqeq100}
 \overline{\mfO} \mfp^{^\mr{Zzz...}}_{\mfs \mfl_2, \mfX} \ \ \ \left(\text{resp.,} \ \mfM \overline{\mfO} \mfp^{^\mr{Zzz...}}_{\mfs \mfl_2, \mfX}\right)
\end{align}
the subset of $ \overline{\mfO} \mfp_{\mfs \mfl_2, \mfX}$ (resp., $\mfM \overline{\mfO} \mfp_{\mfs \mfl_2, \mfX}$) consisting of 
 dormant $\mfs \mfl_2$-opers (resp., dormant generic Miura $\mfs \mfl_2$-opers).
The assignment $(\mcE_B, \nabla_\mcE, \mcE'_B, \eta_\mcE) \mapsto (\mcE_B, \nabla_\mcE)$ determines a map of sets
\begin{align} \label{aaa3}
\mu_\mfX : \mfM \overline{\mfO} \mfp_{\mfs \mfl_2,  \mfX}\migi \overline{\mfO} \mfp_{\mfs \mfl_2, \mfX},
\end{align}
which is nothing but the Miura transformation discussed in Introduction.
This map restricts to a map 
\begin{align} \label{aaa1}
\mu^{^\mr{Zzz...}}_\mfX : \mfM \overline{\mfO} \mfp_{\mfs \mfl_2,  \mfX}^{^\mr{Zzz...}} \migi \overline{\mfO} \mfp_{\mfs \mfl_2, \mfX}^{^\mr{Zzz...}},
\end{align}
referred to as the {\bf dormant Miura transformation} for $\mfX$.

\vspace{5mm}
\subsection{Radius  and exponent}
\leavevmode\\ \vspace{-4mm}

Next, let $\mft$ be the Lie algebra associated with  the maximal torus of $\mr{PGL}_2$ consisting of the image (via the quotient $\mr{GL}_2 \migisurj \mr{PGL}_2$) of diagonal matrices.
Denote by $\mfc$  the GIT quotient $\mfs \mfl_2 \ooalign{$/$ \cr $\,/$}\hspace{-0.5mm}\mr{PGL}_2$ of 
   $\mfs \mfl_2$  by 
the adjoint action of $\mr{PGL}_2$ and 
 by $\chi : \mft \migisurj \mfc$ the composite quotient $\mft \migiincl \mfs \mfl_2 \migisurj \mfc$.
Since $\mft$ and $\mfc$ is defined over $\mbF_p$, 
 it makes sense to speak of the sets $\mft (\mbF_p)$, $\mfc (\mbF_p)$ of the  $\mbF_p$-rational points of $\mft$, $\mfc$ respectively.
Let us write
\begin{align}
\check{\rho} := \begin{pmatrix} \frac{1}{2} & 0 \\ 0 & -\frac{1}{2} \end{pmatrix} 
\in \mft (\mbF_p).
\end{align}
Then, the assignment $\varepsilon \mapsto \varepsilon\cdot \check{\rho}$ determines a bijection $k \isom \mft (k)$.

According to ~\cite{Wak5}, Definition 2.9.1,  the notion of the {\it radii} of a given $\mfs \mfl_2$-oper is defined  as  an element of $\mfc (k)^{\times r}$ (if $r >0$).
Also, the notion of the {\it exponent} (cf. ~\cite{Wak7}, Definition 3.6.1) of a given generic Miura $\mfs \mfl_2$-oper is defined as  an element of  $\mft (k)^{\times r}$ (if $r >0$).
For convenience, we shall refer  to any $\mfs \mfl_2$-oper  (resp., generic Miura $\mfs \mfl_2$-oper) as being  
{\it of radii $\emptyset$} (resp., {\it of exponent $\emptyset$}).
For each $\vec{\rho} \in \mfc (k)^{\times r}$ (resp., $\vec{\varepsilon} \in \mft (k)^{\times r}$), where $\vec{\rho} := \emptyset$ (resp., $\vec{\varepsilon} := \emptyset$) if $r =0$, 
we  denote by 
\begin{align}
 \overline{\mfO} \mfp^{^\mr{Zzz...}}_{\mfs \mfl_2, \mfX, \vec{\rho}} \ \ \ \left(\text{resp.,} \ \mfM \overline{\mfO} \mfp^{^\mr{Zzz...}}_{\mfs \mfl_2, \mfX, \vec{\varepsilon}}\right)
\end{align}
the subset of $ \overline{\mfO} \mfp^{^\mr{Zzz...}}_{\mfs \mfl_2, \mfX}$ (resp., $\mfM \overline{\mfO} \mfp^{^\mr{Zzz...}}_{\mfs \mfl_2, \mfX}$)
consisting of dormant $\mfs \mfl_2$-opers of radius $\vec{\rho}$ (resp., dormant generic Miura $\mfs \mfl_2$-opers of exponent $\vec{\varepsilon}$).
It follows from 
~\cite{Wak5}, Theorem C, and  ~\cite{Wak7}, Theorem 3.8.3, that 
these sets are finite.
Moreover, the sets  $ \overline{\mfO} \mfp^{^\mr{Zzz...}}_{\mfs \mfl_2, \mfX}$ and 
$ \mfM \overline{\mfO} \mfp_{\mfs \mfl_2, \mfX}^{^\mr{Zzz...}}$ decompose into the disjoint unions
\begin{align} \label{EE002}
 \overline{\mfO} \mfp^{^\mr{Zzz...}}_{\mfs \mfl_2, \mfX} = \coprod_{\vec{\rho} \in \mfc (\mbF_p)^{\times r}}  \overline{\mfO} \mfp^{^\mr{Zzz...}}_{\mfs \mfl_2, \mfX, \vec{\rho}}, \hspace{10mm}
 \mfM \overline{\mfO} \mfp_{\mfs \mfl_2, \mfX}^{^\mr{Zzz...}} = \coprod_{\vec{\varepsilon} \in \mft (\mbF_p)^{\times r}} 
\mfM \overline{\mfO} \mfp_{\mfs \mfl_2, \mfX, \vec{\varepsilon}}^{^\mr{Zzz...}}
\end{align}
respectively.
The map (\ref{aaa1}) restricts to a map
\begin{align} \label{aaa2}
\mu^{^\mr{Zzz...}}_{\mfX, \vec{\varepsilon}} : \mfM \overline{\mfO} \mfp_{\mfs \mfl_2,  \mfX, \vec{\varepsilon}}^{^\mr{Zzz...}} \migi \overline{\mfO} \mfp_{\mfs \mfl_2, \mfX, \chi (\vec{\varepsilon})}^{^\mr{Zzz...}},
\end{align}
where $\chi (\vec{\varepsilon}) := (\chi (\varepsilon_i))_{i=1}^r$ if $r>0$ (resp.,  $\chi (\vec{\varepsilon}) :=\emptyset$ if $r =0$).
\vspace{5mm}
\subsection{Pre-Tango structures on a pointed stable curve}
\leavevmode\\ \vspace{-4mm}


In this subsection, we recall the definition of a pre-Tango structure given in ~\cite{Wak7}, Definition 5.3.1.
Let $\mfX := (X, \{ \sigma_i \}_{i=1}^r)$ be as above.

\vspace{3mm}
\bde \label{D170} \leavevmode\\
\ \ \ 
A {\bf pre-Tango structure} on $\mfX$ 
 is defined to be a $k^\mr{log}$-connection $\nabla_\Omega$
 on the line bundle $\Omega_{X^\mr{log}/k^\mr{log}}$ (i.e., a 
$k$-linear morphism $\nabla_\Omega : \Omega_{X^\mr{log}/k^\mr{log}} \migi \Omega_{X^\mr{log}/k^\mr{log}}^{\otimes 2}$ satisfying the Leibniz rule: $\nabla_\Omega (a \cdot  v) = d a \otimes v + a \cdot \nabla_\Omega (v)$, where $a \in \mcO_X$, $v \in \Omega_{X^\mr{log}/k^\mr{log}}$)
satisfying the following two conditions:
\begin{itemize}
\item
It has vanishing $p$-curvature;
\item
If $C_{X^\mr{log}/k^\mr{log}}$ denotes the Cartier operator 
 $\Omega_{X^\mr{log}/k^\mr{log}} \migi \Omega_{X^{(1) \mr{log}}/k^\mr{log}}$
 induced by (the inverse of) ``\,$C^{-1}$\," resulting from  ~\cite{KaKa}, Theorem 4.12 (1),
  where $X^{(1)}$ is the Frobenius twist of $X$ over $k$, then the inclusion relation 
 $\mr{Ker} (\nabla_\Omega) \subseteq \mr{Ker} (C_{X^\mr{log}/k^\mr{log}})$ holds.
\end{itemize}
 \ede
\vspace{3mm}

Note that  it makes sense to speak of the {\it monodromy} (at each marked point $\sigma_i$ of $\mfX$) of a pre-Tango structure on $\mfX$ (cf. ~\cite{Wak5}, Definition 1.6.1 for the definition of monodromy).
The monodromy of a pre-Tango structure lies in 
 $k$ ($\cong \mcE nd_{\mcO_X} (\sigma^*_i (\Omega_{X^\mr{log}/k^\mr{log}}))$).
For convenience, we shall refer to any pre-Tango structure as being  {\it of monodromy $\emptyset$}.

Let $\vec{\varepsilon} := (\varepsilon_i)_{i=1}^r$ be an element of $\mbF_p^{\times r}$ (where  $\vec{\varepsilon} := \emptyset$ if $r =0$).
Denote by
\begin{align}
\overline{\mfT} \mfa \mfn_{\mfX} \ \left(\text{resp.,} \ \overline{\mfT} \mfa \mfn_{\mfX, \vec{\varepsilon}}\right)
\end{align}
the set of pre-Tango strcutres on $\mfX$ (resp., pre-Tango structures on $\mfX$ of monodromy $\vec{\varepsilon}$).
The set $\overline{\mfT} \mfa \mfn_{\mfX}$ decomposes into the disjoint union
\begin{align} \label{TTTT}
\overline{\mfT} \mfa \mfn_{\mfX} = \coprod_{\vec{\varepsilon} \in \mbF_p^{\times r}} \overline{\mfT} \mfa \mfn_{\mfX, \vec{\varepsilon}}
\end{align}
(cf. ~\cite{Wak7}, (190)).
If $\overline{\mfC} \mfo_\mfX$ denotes the set of $k^\mr{log}$-connections on $\Omega_{X^\mr{log}/k^\mr{log}}$, then 
there exists a natural bijection of sets
\begin{align} \label{GGGG}
\overline{\mfC} \mfo_\mfX \isom \mfM \overline{\mfO} \mfp_{\mfs \mfl_2, \mfX},
\end{align}
which induces, by restriction, 
 a bijection 
\begin{align} \label{EErt6}
\overline{\mfT} \mfa \mfn_{\mfX} \isom \mfM \overline{\mfO} \mfp_{\mfs \mfl_2, \mfX}^{^\mr{Zzz...}}
\end{align}
(cf. ~\cite{Wak7}, Theorems 4.4.1 and  5.4.1).
If we write 
\begin{align} \label{Ww400}
\vec{\varepsilon} \cdot \check{\rho} := (\varepsilon_i \cdot \check{\rho})_{i=1}^r
\end{align}
 (where $\vec{\varepsilon} \cdot \check{\rho} := \emptyset$ if $r=0$), 
then  it restricts to a bijection 
\begin{align} \label{EE007}
\overline{\mfT} \mfa \mfn_{\mfX, \vec{\varepsilon}} \isom \mfM \overline{\mfO} \mfp_{\mfs \mfl_2, \mfX, \vec{\varepsilon}\cdot \check{\rho}}^{^\mr{Zzz...}}.
\end{align}

\vspace{5mm}
\subsection{Gluing pre-Tango structures}
\leavevmode\\ \vspace{-4mm}

Let  us discuss the procedure for gluing  pre-Tango structures by means of a clutching data.
To begin with,  we shall  define the notion of  a clutching data,  as follows.

\bde \leavevmode\\
 \ \ \ 
 A  {\bf clutching  data  of type $(g,r)$}  is   a collection of data:
\begin{equation}  \label{Er4}
 \mcD :=  ( \Gamma^+, \{ (g_v, r_v) \}_{v \in V_\Gamma}, \{ \lambda_v \}_{v \in V_\Gamma}),  
  \end{equation}
where

\begin{itemize}
\item[$\bullet$]
$\Gamma^+ := (V_\Gamma, E_\Gamma, \zeta_\Gamma, \lambda_\Gamma)$ denotes a marked semi-graph with   $r = \sharp (B_{V_\Gamma})$ whose underlying semi-graph is    (finite and) connected;  
\item[$\bullet$]
$(g_v, r_v)$ (for each $v  \in V_\Gamma$) denotes   a pair of nonnegative integers with $2g_v -2 + r_v >0$, $r_v >0$, 
 and $g = \beta (\Gamma) + \sum_{j =1}^n g_j$;
\item[$\bullet$]
$\lambda_v$ (for each $v  \in V_\Gamma$) denotes  a bijection  $\lambda_v : B_{v} \isom  \{ 1, \cdots, r_v \}$ of sets.
\end{itemize} 
\vspace{3mm}
\ede

 
Let $\mcD$ be a clutching data  of type $(g,r)$ as of   (\ref{Er4}) and $\{ {\mfX_v} \}_{v \in V_\Gamma}$
 a collection of pointed stable curves   over $k$ indexed by the elements of $V_\Gamma$, where  each  ${\mfX_v} := (X_v, \{ \sigma_{v, i}  \}_{i=1}^{r_v})$ is  of type $(g_v,r_v)$.
Here,  we  assume that  the curves $X_v$ are all {\it smooth}.
The  pointed  curves ${\mfX_v}$ 
 may be glued together, by means of $\mcD$,  to a new pointed stable curve $\mfX := (X, \{\sigma_i  \}_{i=1}^r)$ of type $(g,r)$ 
in such a way that
\begin{itemize}
\item
the dual marked semi-graph associated with $\mfX$ is given by $\Gamma^+$, where each  vertex $v \in V_\Gamma$ corresponds to the  irreducible component  $X_v$;
\item
if $e = \{ b_1, b_2 \}$ is an edge of $\Gamma$ with $\zeta_\Gamma (b_1) = v_{1}$, $\zeta_\Gamma (b_2) =v_{2}$ (for some $v_1$, $v_2 \in V_\Gamma$), then $e$ corresponds to the  node of $X$ obtained by gluing together $X_{v_1}$ at the $\lambda_{v_1} (b_1)$-th marked point $\sigma_{v_1, \lambda_{v_1} (b_1)}$  to $X_{v_2}$ at the $\lambda_{v_2} (b_2)$-th marked point $\sigma_{v_2, \lambda_{v_2} (b_2)}$;
\item
  the $i$-th ($i = 1, \cdots, r$)   marked point of $\mfX$   arises from the $\lambda_{\zeta_\Gamma (\lambda_{\Gamma}^{-1} (i)^\divideontimes)} (\lambda_{\Gamma}^{-1} (i)^\divideontimes)$-th marked point of $\mfX_{\zeta_\Gamma (\lambda_{\Gamma}^{-1} (i)^\divideontimes)}$.
\end{itemize}
One may extend, in the evident way,  this construction to the case where  $X_v$'s are (possibly {\it non-smooth}) pointed stable curves.

Denote by $\mfC \mfl \mfu \mft_v$   ($v \in V_\Gamma$) the resulting morphism  $X_v \migi X$.
We shall write  $X_v^{\mfX \text{-} \mr{log}}$ for the  log scheme obtained by equipping $X_v$ with the log structure pulled-back from the log structure of $X^\mr{log}$ via $\mfC \mfl \mfu \mft_v$.
The structure morphism $X_v \migi \mr{Spec} (k)$ of  $X_v$ extends to a morphism $X_v^{\mfX \text{-} \mr{log}} \migi \mr{Spec} (k)^{\mfX \text{-}\mr{log}}$ of log schemes.
Moreover, the morphism $X_v^{\mfX \text{-} \mr{log}} \migi X_v$  
extends naturally to a commutative square diagram
\begin{align}
\begin{CD}
X_v^{\mfX \text{-} \mr{log}}@>>> X_v^{\mfX_v \text{-} \mr{log}}
\\
@VVV @VVV
\\
\mr{Spec} (k)^{\mfX \text{-} \mr{log}}@>>> \mr{Spec} (k)^{\mfX_v \text{-} \mr{log}}.
\end{CD}
\end{align}
The resulting morphism
\begin{align}
\mfe_v : X_v^{\mfX \text{-} \mr{log}} \migi X_v^{\mfX_v \text{-} \mr{log}} \times_{\mr{Spec} (k)^{\mfX_v \text{-}\mr{log}}} \mr{Spec} (k)^{\mfX \text{-}\mr{log}}
\end{align}
 is verified to be  
 {\it log \'{e}tale}.
(Note that the underlying morphism  of $\mfe_v$ coincides with the identity morphism of $X_v$.)
In particular, it induces an isomorphism
 \begin{align} \label{E30}
 \Omega_{X_{v}^{\mfX_v \text{-}\mr{log}}/k^{\mfX_v \text{-}\mr{log}}} \isom \Omega_{X_v^{\mfX \text{-} \mr{log}}/k^{\mfX \text{-} \mr{log}}} \left(\isom \mfC \mfl \mfu \mft_{v}^*(\Omega_{X^{\mfX \text{-}\mr{log}}/k^{\mfX \text{-}\mr{log}}})\right)
 \end{align}
  of $\mcO_{X_v}$-modules.

Next,  let  $\nabla_\Omega$
be a $k^{\mfX \text{-} \mr{log}}$-connection on $\Omega_{X^{\mfX \text{-} \mr{log}}/k^{\mfX \text{-} \mr{log}}}$.
For each $v\in V_\Gamma$, the pull-back   $\mfC \mfl \mfu \mft^*_v (\nabla_\Omega)$
 of $\nabla_\Omega$ to $X_v$ forms  a $k^{\mfX \text{-} \mr{log}}$-connection  on $\Omega_{X_v^{\mfX \text{-} \mr{log}}/k^{\mfX \text{-} \mr{log}}}$.
 Moreover, 
 by means of  (\ref{E30}),
  $\mfC \mfl \mfu \mft^*_v (\nabla_\Omega)$ may be thought of as  a $k^{\mfX_v \text{-} \mr{log}}$-connection 
  \begin{align} \label{fff05}
  \nabla_\Omega^v :  \Omega_{X_{v}^{\mfX_v \text{-}\mr{log}}/k^{\mfX_v\text{-}\mr{log}}} \migi  \Omega^{\otimes 2}_{X_{v}^{\mfX_v \text{-}\mr{log}}/k^{\mfX_v \text{-}\mr{log}}}
  \end{align}
   on $\Omega_{X_{v}^{\mfX_v \text{-}\mr{log}}/k^{\mfX_v \text{-}\mr{log}}}$.
    We shall refer to $\nabla_\Omega^v$ as the {\bf restriction} of $\nabla_\Omega$ to $\mfX_v$.

\bpr \label{P04gg} \leavevmode\\
 \ \ \ 
 Let us keep the above  notation.
 \begin{itemize}
 \item[(i)]
Suppose  that $\nabla_\Omega$ specifies   a pre-Tango structure on $\mfX$.
Then,   for each $v \in V_\Gamma$, the restriction  $\nabla_{\Omega}^v$
   specifies  a pre-Tango structure  on $\mfX_v$.
If, moreover,  $\nabla^v_\Omega$ ($v \in V_\Gamma$) is of monodromy  $\vec{\varepsilon}^{\,v} := (\varepsilon_i^v)_{i=1}^{r_v}\in \mbF_p^{\times r_v}$, 
then the collection  
\begin{align}
\left(\tau^{-1} (\varepsilon^{\zeta_\Gamma (b)}_{\lambda_{\zeta_\Gamma (b)} (b)}) \right)_{b \in B_\Gamma} \in \widetilde{\mbF}_p^{B_\Gamma}
\end{align}
 (where $\varepsilon^{V_\Gamma}_{\lambda_{V_\Gamma}(b)} := -\varepsilon^{\zeta_\Gamma (b^\divideontimes)}_{\lambda_{\zeta_\Gamma (b^\divideontimes)} (b^\divideontimes)}$ for any $b \in B_{V_\Gamma}$)
  forms  a $p$-branch numbering on $\Gamma$.
 \item[(ii)]
 Conversely, let  $ \vec{m}_\mcD := (m_b )_{b \in B_\Gamma} \in \widetilde{\mbF}_p^{B_\Gamma}$ be  a $p$-branch numbering on $\Gamma$ of exponent $\vec{\varepsilon} \in \mbF_p^{\times r}$ (where $\vec{\varepsilon} := \emptyset$ if $r =0$).
 Suppose that  on each ${\mfX_v}$, 
we are given a pre-Tango structure $\nabla_{\Omega, v}$  of monodromy 
$(\tau (m_{\lambda^{-1}_v(i)}))_{i=1}^{r_v}$.
Then, there exists a pre-Tango structure $\nabla_\Omega$ on $\mfX$
 of monodromy  $\vec{\varepsilon}$,
 which is uniquely determined   by the condition that  for any $v \in V_\Gamma$ the  restriction 
of $\nabla_\Omega$ to ${\mfX_v}$
 coincides with    $\nabla_{\Omega, v}$.

  \end{itemize}

 \epr
\begin{proof}
Assertions (i) and (ii) follow immediately  from ~\cite{Wak5}, Proposition 7.2.1 (i), (ii) and the fact that for each $v \in V_\Gamma$, 
the restriction of the Cartier operator $C_{X^{\mfX \text{-} \mr{log}}/k^{\mfX \text{-} \mr{log}}}$ to $X_j$ may be identified, via  (\ref{E30}), with $C_{X_v^{\mfX_v \text{-} \mr{log}}/k^{\mfX_v \text{-} \mr{log}}}$.
\end{proof}
\vspace{3mm}

Let $\vec{\varepsilon}$ be an element of $\mbF_p^{\times r}$  (where  $\vec{\varepsilon} := \emptyset$ if $r =0$),
$\mcD$  a clutching data of type $(g,r)$ with  underlying marked semi-graph  $\Gamma^+ := (V_\Gamma, E_\Gamma, \zeta_\Gamma, \lambda_\Gamma)$,
 and 
 $\vec{m}_\mcD := (m_b)_{b \in B_\Gamma} \in \widetilde{\mbF}_p^{B_\Gamma}$  a $p$-branch numbering on $\Gamma$ of exponent $\vec{\varepsilon}$.
If we write
$\vec{\varepsilon}_\mcD^{\ v} := (\tau (m_{\lambda^{-1}_v(i)}))_{i=1}^{r_v}$ (for each $v \in V_\Gamma$), then
 (by applying Proposition \ref{P04gg} above) a collection of pre-Tango structures $(\nabla_{\Omega, v})_{v \in V_\Gamma}$, where 
each $\nabla_{\Omega, v}$ is of monodromy $\vec{\varepsilon}_\mcD^{\ v}$,
induces a pre-Tango structure $\nabla_\Omega$ on $\mfX$ of monodromy $\vec{\varepsilon}$.
 The assignment  $(\nabla_{\Omega, v})_{v \in V_\Gamma} \mapsto \nabla_\Omega$ 
 induces 
 a morphism
 \begin{align}
 \mfC \mfl \mfu \mft_{\mcD, \vec{m}_\mcD} : 
 \prod_{v \in V_\Gamma} \overline{\mfT} \mfa \mfn_{\mfX_v, \vec{\varepsilon}_\mcD^{\, v}} \migi \overline{\mfT} \mfa \mfn_{\mfX, \vec{\varepsilon}}.
 \end{align}

\bpr \label{yy0176}\leavevmode\\
 \ \ \ 
Let  $\mcD$ and $\vec{\varepsilon}$ be as above.
Then,  the following map is bijective:
\begin{align}
\coprod_{\vec{m}_\mcD} \prod_{v \in V_\Gamma} \overline{\mfT} \mfa \mfn_{\mfX_v, \vec{\varepsilon}_\mcD^v} \xrightarrow{\coprod \mfC \mfl \mfu \mft_{\mcD, \vec{m}_\mcD}} \overline{\mfT} \mfa \mfn_{\mfX, \vec{\varepsilon}},
\end{align}
 where the disjoint union 
 in the left-hand side
  is taken over the  set of 
 $p$-branch numberings  on $\Gamma$ of exponent $\vec{\varepsilon}$.
 \epr
\begin{proof}
The assertion follows immediately from the decomposition (\ref{EE002}), Proposition \ref{P04gg},  and the definition of $\mfC \mfl \mfu \mft_{\mcD, \vec{m}_\mcD}$.
\end{proof}

\vspace{10mm}
\section{Combinatorial  description of dormant Miura opers} \vspace{3mm}



\vspace{5mm}
\subsection{Balanced $p$-edge numberings} \label{SeSe335}
\leavevmode\\ \vspace{-4mm}

In this section, we shall study a 
combinatorial
 description of dormant Miura $\mfs \mfl_2$-opers  (equivalently, pre-Tango structures) on a  totally degenerate curve.
 We first recall the 
 combinatorial
  description of dormant $\mfs \mfl_2$-opers on a totally degenerate curve, which
 was essentially given in the previous work of  $p$-adic Teichm\"{u}ller theory due to S. Mochizuki (cf. ~\cite{Mzk2}, Chap.\,V, \S\,1, (3), p.\,232, or ~\cite{Wak5},  \S\,7.11). 
The  objects  used to describe  combinatorially  dormant $\mfs \mfl_2$-opers  will be referred, in the present paper,  to as {\it balanced $p$-edge numberings} (cf. Definition \ref{De113} below).
Let us fix 
a marked semi-graph $\Gamma^+ := (V_\Gamma, E_\Gamma, \zeta_\Gamma, \lambda_\Gamma)$
whose underlying semi-graph is   connected, $3$-regular, and of type $(g,r)$.

\bde \label{De113}\leavevmode\\
\vspace{-5mm}
\begin{itemize}
\item[(i)]
 A {\bf balanced $p$-edge numbering} on $\Gamma^+$ is 
 a collection $\vec{m} := (m_b)_{b \in B_\Gamma} \in \widetilde{\mbF}_p^{B_\Gamma}$ of elements of $\widetilde{\mbF}_p$ indexed by $B_\Gamma$ 
 satisfying the following two conditions: 
 \begin{itemize}
 \item[$\bullet$]
For each edge $e := \{ b, b^\divideontimes\} \in E_{\Gamma}$, the equality $m_b = m_{b^\divideontimes}$ holds.
 \item[$\bullet$]
 For each vertex $v \in V_\Gamma$ (where we write  $B_v := \{b_1, b_2, b_3 \}$ and $m_l := m_{b_l}$ for each  $l =1,2,3$),
   the inequalities in   $\bigstar_{m_1, m_2, m_3}$
   displayed below
     are satisfied:
\vspace{2mm}
\begin{itemize}
\item[]
 \begin{itemize}
 \item[$\bigstar_{m_1, m_2, m_3}$ : ]
$|m_2 -m_3| \leq m_1 \leq m_2 +m_3$,
 \hspace{5mm} $m_1 + m_2 + m_3 \leq  p-2$.
\end{itemize}
 \end{itemize}
\vspace{2mm}
 \end{itemize}
(By the first condition, any balanced $p$-edge numbering may be thought of as a  numbering described  on the set of edges $E_\Gamma$, as its name suggests.)
 \item[(ii)]
 Assume further  that $B_{V_\Gamma} \neq \emptyset$.
 Let $\vec{m} := (m_b)_{b \in B_{V_\Gamma}} \in \widetilde{\mbF}_p^{B_\Gamma}$ be a balanced $p$-edge numbering on $\Gamma$ and $\vec{\varepsilon} := (\varepsilon_i)_{i=1}^r$ an element of $\mbF_p^{\times r}$.
 We shall say that $\vec{m}$ is {\bf of radii  $\vec{\varepsilon}$} if $\tau (m_b) = \varepsilon_{\lambda_\Gamma (b)}$ for any $b \in B_{V_\Gamma}$. For convenience, (regardless of whether $B_{V_\Gamma}$ is empty or not) we shall refer to any $p$-edge numbering as being {\bf of radii $\emptyset$}. 
 \end{itemize}
 \ede
\vspace{3mm}

Denote by 
\begin{align}
p \text{-}\mfE \mfd_{\Gamma^+}^\mr{bal}
\end{align}
the set of balanced $p$-branch numberings  on $\Gamma^+$.
Also, for each  $\vec{\varepsilon} \in \mbF_p^{\times r}$ (where  $\vec{\varepsilon} := \emptyset$ if $r =0$), 
 we shall write
\begin{align} \label{eq5567}
p \text{-} \mfE \mfd_{\Gamma^+, \vec{\varepsilon}}^\mr{bal}
\end{align}
for the subset of  $p \text{-}\mfE \mfd_{\Gamma^+}^\mr{bal}$ consisting of balanced $p$-branch numberings of radii  $\vec{\varepsilon}$.
 In particular, we have 
 \begin{align} \label{YYYY}
 p \text{-}\mfE \mfd_{\Gamma^+}^\mr{bal} = \coprod_{\vec{\varepsilon} \in \mbF_p^{\times r}} p \text{-} \mfE \mfd_{\Gamma^+, \vec{\varepsilon}}^\mr{bal}.
 \end{align}

\vspace{5mm}
\subsection{Combinatorial description of  dormant $\mfs \mfl_2$-opers} \label{sc3345}
\leavevmode\\ \vspace{-4mm}

Next, let us construct a bijective correspondence between the  set of dormant $\mfs \mfl_2$-opers on a  totally degenerate curve  
 and the set of balanced $p$-branch numberings  on the dual semi-graph of this curve.
First, let us consider the case where $(g,r) = (0, 3)$.
Denote by $[0]$, $[1]$, and $[\infty]$ the $k$-rational points of  the projective line $\mbP^1$  over $k$ determined by the values $0$, $1$, and $\infty$ respectively.
After ordering  the points $ [0], [1], [\infty]$ suitably (say, $\sigma_1 := [0]$, $\sigma_2 := [1]$, and $\sigma_3 := [\infty]$), we obtain 
 a unique (up to isomorphism) pointed stable curve 
\begin{equation} \label{1051}
\mfP := (\mbP^1, \{ \sigma_i\}_{i=1}^3)
\end{equation}
 of type $(0,3)$ over $k$.
The dual marked semi-graph $\Gamma^+_\mfP := (V_\mfP, E_\mfP, \zeta_\mfP, \lambda_\mfP)$ associated with $\mfP$ is given as follows:
\begin{itemize}
\item[$\bullet$]
$V_\mfP := \left\{ v \ (:= \mbP^1) \right\}$;
\item[$\bullet$]
$E_\mfP := \{ e_1, e_2, e_3 \}$, where 
\begin{align}
e_1 := \{ \sigma_1, \{ \sigma_1 \} \}, \ \ e_2 := \{ \sigma_2, \{ \sigma_2 \} \}, \ \ e_3 := \{ \sigma_3, \{ \sigma_3 \} \};
\end{align}
\item[$\bullet$]
$\zeta_\mfP : B_{\Gamma_\mfP} \migi  V_\mfP \sqcup \{ V_\mfP \}$ is given by 
\begin{align}
\zeta_\mfP (\sigma_1) = \zeta_\mfP ( \sigma_2) = \zeta_\mfP ( \sigma_3) = v,  \ \ \ \zeta_\mfP ( \{ \sigma_1 \}) = \zeta_\mfP ( \{ \sigma_2 \}) = \zeta_\mfP ( \{ \sigma_3 \}) = V_\mfP.
\end{align}
\item[$\bullet$]
$\lambda_\mfP : B_{V_\mfP} \ \left(= \{  \{ \sigma_1 \},  \{ \sigma_2 \},  \{ \sigma_3 \} \}\right) \isom \{ 1, 2,3\}$ is given by $ \{ \sigma_l \} \mapsto l$ for any $l \in \{ 1,2,3\}$.
\end{itemize}
The assignment $(m_{b})_{b \in B_{\Gamma_\mfP}} \mapsto (m_{\sigma_1}, m_{\sigma_2}, m_{\sigma_3})$ 
 gives a bijective correspondence between   the set of  balanced $p$-edge numberings on $\Gamma_\mfP^+$ and   the set of  triples $(m_1, m_2, m_3) \in \mbZ^{\times 3}$ of  integers satisfying 
the condition  $\bigstar_{m_1, m_2, m_3}$.
The inverse assignment is given by $(m_1, m_2, m_3) \mapsto (m_{b})_{b \in B_\Gamma}$, where $m_{\sigma_l} := m_l$, $m_{ \{ \sigma_l \}} := m_l$
  (for any $l \in \{1, 2,3\}$).
By passing to this correspondence,  {\it we shall identify   each balanced $p$-edge numbering on $\Gamma_\mfP^+$ with  such a triple $(m_1, m_2, m_3)$}.

Now, 
denote by $\iota$ the bijection defined as 
\begin{align}
\hspace{20mm} \iota : 
k & \migi \hspace{10mm} \mfc (k) \\
\vin 
 & 
 \hspace{19mm}
 \vin \notag \\
a  
 & \mapsto 
 \chi \left( \begin{pmatrix} a +  \overline{\frac{p+1}{2}} & 0 \\ 0 &- \left( a + \overline{\frac{p+1}{2}}\right)\end{pmatrix}\right), \notag 
\end{align}
 where $\overline{(-)}$ denotes the image of $(-)$ via the composite  $\mbZ \migisurj \mbF_p \migiincl k$. 
Recall from    ~\cite{Mzk2}, Chap.\,I, \S\,4.3, p.\,117, Theorem 4.4 (or ~\cite{Wak5}, Theorem A), that  for each $(\rho_1, \rho_2, \rho_3) \in \mfc (k)^{\times 3}$
there exists a unique $\mfs \mfl_2$-oper $\mcE^\spadesuit_{\rho_1, \rho_2, \rho_3}$ on $\mfP$ of radii $(\rho_1, \rho_2, \rho_3)$.
That is to say, 
we obtain a bijection
\begin{align} \label{aaa6}
k^{\times 3} \hspace{3mm} &\isom \overline{\mfO} \mfp_{\mfs \mfl_2, \mfP} \\
\vin \hspace{7mm} & \hspace{8mm} \vin \notag \\
(a_1, a_2, a_3) &\mapsto \mcE^\spadesuit_{\iota (a_1), \iota (a_2), \iota (a_3)}. \notag
\end{align}

\ble \label{LLL001}\leavevmode\\
 \ \ \
The composite $\mbZ^{\times 3} \migisurj \mbF_p^{\times 3}\migiincl  k^{\times 3} \xrightarrow{(\ref{aaa6})} \overline{\mfO} \mfp_{\mfs \mfl_2, \mfP}$ restricts to a bijection
\begin{align} \label{eq2233}
p\text{-} \mfE \mfd_{\Gamma^+_\mfP}^\mr{bal} \isom \overline{\mfO} \mfp_{\mfs \mfl_2,  \mfP}^{^\text{Zzz...}}.
\end{align}
 \ele
\begin{proof}
The assertion follows from   ~\cite{Mzk2}, Chap.\,V, \S\,1, (3), p.\,232 (cf. ~\cite{Wak5}, the discussion in \S\,7.11.
\end{proof}
\vspace{3mm}

Next, we shall extend  the above   result to the case where the underlying curve is 
an arbitrary totally degenerate curve.
To this end, let us recall the definition of a  totally degenerate curve.
Let  $\mfX : = (X, \{ \sigma_i \}_{i =1}^r)$ be  a pointed stable curve over $k$ of type $(g,r)$.
Write 
 $\nu_{\mfX} : \coprod_{l =1}^{L_{\mfX}} X_l \migi X$ for  the normalization of $X$, where  $L_{\mfX}$ denotes some positive integer and  each $X_l$ ($l = 1, \cdots, L_{\mfX_{/k}}$) is a proper {\it smooth} connected curve over $k$.
Then, we shall say that  $\mfX$ is {\bf totally degenerate} if, for any $l =1, \cdots, L_{\mfX}$, the pointed stable curve
\begin{equation}
{\mfX_l} := (X_l, \nu_{\mfX}^{-1}(E_{\mfX}) \cap X_l (k))
\end{equation}
 is isomorphic to $\mfP$, where  we consider  $E_{\mfX} := N_\mfX \sqcup \{ \sigma_i \}_{i=1}^r$ (cf. (\ref{EEE445}))   as a subset of $X (k)$.

Now, let $\mfX$ be a  totally degenerate pointed stable curve over $k$ of type $(g,r)$.
Then, $\mfX$ may be obtained by gluing together   finite copies  of $\mfP$ by means of some clutching data whose underlying  marked semi-graph is $\Gamma_\mfX^+$.
For each $v \in V_\mfX$, we shall denote by $\mfP_v$ the $3$-pointed projective line corresponding to $v$ (i.e., $\mfP_v \cong \mfP$).
Also, for each $v \in V_\mfX$ and $b \in B_v$, we shall denote by $\sigma_b$ the marked point of $\mfP_v$ corresponding to $b$.

Let $(m_b)_{b \in B_{\Gamma_\mfX}}$ be a balanced $p$-edge numbering on $\Gamma_\mfX^+$. 
For each  $v \in V_\mfX$ with $B_v := \{ b_1, b_2, b_3\}$, the triple $(m_{b_1}, m_{b_2}, m_{b_3})$ specifies a  balanced $p$-edge numbering on $\Gamma_{\mfP_v}^+$.
This triple corresponds, via   (\ref{eq2233}),  to a dormant $\mfs \mfl_2$-oper  $\mcE^\spadesuit_v$ on $\mfP_v$.
One may assume, without loss of generality, that each $\mcE_v^\spadesuit$ is of canonical type (cf. ~\cite{Wak5}, Definition 2.7.1).
If $e := \{ b, b^\divideontimes \} \in E_{\Gamma_\mfX}$ is an edge with $\{ V_\Gamma \}  \notin \zeta_\Gamma (e)$,
then the radius    of $\mcE^\spadesuit_{\zeta_{\Gamma_\mfX}  (b)}$ at $\sigma_{b}$ coincides with  the radius 
of  $\mcE^\spadesuit_{\zeta_{\Gamma_\mfX}(b^\divideontimes)}$ at $\sigma_{b^\divideontimes}$.
It  follows from ~\cite{Wak5}, Proposition 7.3.3 (ii) that $\mcE^\spadesuit_v$'s may be glued together to a dormant $\mfs \mfl_2$-oper $\mcE^\spadesuit$ on $\mfX$ (of canonical type).
 The bijectivity of  (\ref{eq2233}) implies  the following proposition.

\vspace{3mm}
\bpr \label{y0176}\leavevmode\\
 \ \ \ 
Let $\mfX$ be as above.
Then, the assignment 
$(m_b)_{b \in B_{\Gamma_\mfX}} \mapsto \mcE^\spadesuit$
 discussed above defines a bijection
\begin{align} \label{eq889}
 p \text{-}\mfE \mfd^\mr{bal}_{\Gamma^+_\mfX} \isom \mfO \mfp_{\mfs \mfl_2, \mfX}^{^\mr{Zzz...}}.
\end{align}
If, moreover, $r >0$, then for each 
$\vec{\varepsilon} := (\varepsilon_i)_{i=1}^r\in \mbF_p^{\times r}$, 
the bijection (\ref{eq889}) restricts to a bijection
\begin{align} \label{eq890}
p \text{-}\mfE \mfd^\mr{bal}_{\Gamma^+_\mfX, \vec{\varepsilon}} \isom \mfO \mfp_{\mfs \mfl_2, \mfX, \iota (\vec{\varepsilon})}^{^\mr{Zzz...}},
\end{align}
where $\iota (\vec{\varepsilon}) := (\iota (\varepsilon_i))_{i=1}^r$.
 \epr

\vspace{3mm}
\begin{rema} \label{R0ff011}
\leavevmode\\
 \ \ \ 
If we use the notation $(-)^{\, \mu}$  defined in (\ref{eq4729}) below,
then  (\ref{eq890}) may be expressed as
\begin{align} \label{eq8901}
p \text{-}\mfE \mfd^\mr{bal}_{\Gamma^+_\mfX,-\vec{\varepsilon}^{\, \mu}} \isom \mfO \mfp_{\mfs \mfl_2, \mfX, \chi (\vec{\varepsilon}\cdot \check{\rho})}^{^\mr{Zzz...}}
\end{align}
for each  $\vec{\varepsilon} \in \mbF_p^{\times r}$.
 \end{rema}

\vspace{5mm}
\subsection{Strict $p$-branch numberings} \label{sc778}
\leavevmode\\ \vspace{-4mm}

Next, in order to describe  combinatorially dormant generic Miura $\mfs \mfl_2$-opers (equivalently, pre-Tango structures),
  we shall introduce the notion of a strict $p$-branch numbering, as follows.
Let $\Gamma^+ : = (V_\Gamma, E_\Gamma, \zeta_\Gamma, \lambda_\Gamma)$ be a marked semi-graph whose underlying semi-graph $\Gamma$ is   connected,  $3$-regular, and of type $(g,r)$.

\vspace{3mm}
\bde \label{FDE2}\leavevmode\\
 \ \  \ A   {\bf  strict $p$-branch numbering} on $\Gamma^+$ is 
 a $p$-branch numbering  $\vec{m} := (m_b)_{b \in B_\Gamma} \in \widetilde{\mbF}_p^{B_{\Gamma}}$
 on $\Gamma$ with $m_b \neq 0$ (for any $b \in B_\Gamma$) such that  for each vertex $v \in V_\Gamma$  (where we shall write $B_v := \{ b_1, b_2, b_3\}$), the equality
 \begin{align}
 \sum_{j=1}^3 m_{b_j} = 1 +p
 \end{align}
  holds.
 \ede
\vspace{3mm}

Denote by 
\begin{align}  
  p\text{-}\mfB  \mfr_{\Gamma^+}^\mr{st}
\end{align}
 the set of  strict    $p$-branch numberings on $\Gamma^+$.
Also, for each $\vec{\varepsilon} \in \mbF_p^{\times r}$ (where $\vec{\varepsilon} := \emptyset$ if $r =0$),  
 we shall write
\begin{align}
   p\text{-}\mfB  \mfr_{\Gamma^+, \vec{\varepsilon}}^\mr{st}
\end{align}
for the set of  strict $p$-branch numberings on $\Gamma^+$ of exponent  $\vec{\varepsilon}$.
The set  $p\text{-}\mfB  \mfr_{\Gamma^+}^\mr{st}$ decomposes into the disjoint union
\begin{align} \label{UUUU}
p\text{-}\mfB  \mfr^\mr{st}_{\Gamma^+} = \coprod_{\vec{\varepsilon} \in \mbF_p^{\times r}}  p\text{-}\mfB  \mfr^\mr{st}_{\Gamma^+, \vec{\epsilon}}.
\end{align}


Now, let us construct an assignment from each strict  $p$-branch numbering  to a balanced  $p$-edge numbering.
Given an element $m \in \widetilde{\mbF}_p$,  we shall write $m^{\mu}$ for the element of $\widetilde{\mbF}_p$ defined as follows:
\begin{align} \label{eq4729}
m^{\mu} :=  \begin{cases} \frac{p-m -1}{2} & \text{if $m$ is even;} \\
\frac{m-1}{2}  & \text{if $m$ is odd}.
\end{cases}
\end{align}
For each 
 strict $p$-branch numbering   $\vec{m} := (m_b)_{b \in B_\Gamma}$ on $\Gamma^+$,
 the collection $\vec{m}^{\mu} : = (m^{\mu}_b)_{b \in B_\Gamma}$ is verify to  specify  a balanced $p$-edge numbering   on $\Gamma^+$.
 Thus, we obtain  a map of sets
 \begin{align} \label{eq4730}
\mu^\mr{comb}_{\Gamma^+} : p\text{-}\mfB  \mfr_{\Gamma^+}^\mr{st} & \migi p\text{-}\mfE \mfd^\mr{bal}_{\Gamma^+} \\
\vin \  \, & \hspace{10mm}\vin \notag \\
 \vec{m} \ & \mapsto \hspace{3mm} \vec{m}^{\mu}, \notag 
 \end{align}
which we shall refer to as the {\bf combinatorial dormant Miura transformation} for $\Gamma$.
Given  $\vec{\varepsilon} := (\varepsilon_i)_{i=1}^r \in \mbF_p^{\times r}$ (where  $\vec{\varepsilon} : = \emptyset$ if $r =0$), 
the map $\mu^\mr{comb}_{\Gamma} $ restricts to a map
\begin{align}
\mu^\mr{comb}_{\Gamma^+, \vec{\varepsilon}}  : p\text{-}\mfB  \mfr_{\Gamma^+, \vec{\varepsilon}}^\mr{st}  \migi p\text{-}\mfE \mfd_{\Gamma^+, \vec{\varepsilon}^{\, \mu}}^\mr{bal},
\end{align}
where $\vec{\varepsilon}^{\, \mu} := (\varepsilon_i^{\mu})_{i=1}^r$ (and $\vec{\varepsilon}^{\, \mu} := \emptyset$ if $r =0$).

\vspace{5mm}
\subsection{Combinatorial description of dormant generic Miura $\mfs \mfl_2$-opers}
\leavevmode\\ \vspace{-4mm}

We shall describe the relationship between the set of pre-Tango structures and the set of strict $p$-branch numberings.
In a fashion analogous to  the case of balanced $p$-edge numberings, {\it we shall identify  each   $p$-branch numbering $(m_b)_{b \in B_{\Gamma_\mfP}}$ on $\Gamma_\mfP$ with  a triple $(m_1, m_2, m_3) \in \widetilde{\mbF}_p^{\times 3}$, where $m_l := m_{\sigma_l} = m^\veebar_{\{ \sigma_l \}}$ ($l =1,2,3$)}.

\vspace{3mm}
\ble \label{PP004} \leavevmode\\
 \ \ \ 
 Let us consider the set
 \begin{align}
p\text{-}\mfB \mfr_{\Gamma^+_\mfP} :=  \left\{(m_1, m_2, m_3) \in \widetilde{\mbF}_p^{\times 3} \ \Bigg| \ p  | \left(\sum_{i=1}^3 m_i -1\right) \right\}.
 \end{align}
Then, an element $(m_1, m_2, m_3)$ of $p\text{-}\mfB \mfr_{\Gamma^+_\mfP}$ specifies a strict $p$-branch numbering on $\Gamma^+_\mfP$ (in the above sense) if and only if the triple  $(m_1^{\mu}, m^{\mu}_2, m^{\mu}_3)$ specifies a balanced $p$-edge numbering on $\Gamma^+_\mfP$ (i.e., satisfies the inequalities in $\bigstar_{m_1^{\mu}, m^{\mu}_2, m^{\mu}_3}$).
\ele
\begin{proof}
Let $(m_1, m_2, m_3)$ be   an element of $p\text{-}\mfB \mfr_{\Gamma^+_\mfP}$ such that $(m_1^{\mu}, m^{\mu}_2, m^{\mu}_3)$  satisfies  the inequalities in $\bigstar_{m_1^{\mu}, m^{\mu}_2, m^{\mu}_3}$.
To complete the proof, 
it suffices to verify  that
this  triple specifies  a strict $p$-branch numbering on $\Gamma^+_\mfP$.
Define $m$ to be the integer with $\sum_{i=1}^3 m_i = mp +1$.
Since $m_i \leq p-1$, the inequality $m < 3$ holds.
Consider the case where $m =0$.
Then, after possibly change of ordering, we may assume that $m_1= m_2 =0$ and $m_3 =1$. 
Then, 
\begin{align}
p-2 \geq \sum_{i=1}^3 m_i^{\mu} = \frac{p-1}{2} + \frac{p-1}{2} + 0 = p-1,
\end{align}
 which is a contradiction.
Next, consider the case where $m =2$.
As the sum $\sum_{i=1}^3 m_i$ is odd, either one of the following two cases  (a), (b) is satisfied:
(a) the three integers $m_1$, $m_2$, $m_3$ are all  odd; (b) 
 two of the three integers $m_1$, $m_2$, $m_3$ are even and the remaining one is odd.
But, in  the case (a), we obtain  a contradiction since 
\begin{align}
p-2 \geq  \sum_{i=1}^3 m_i^{\mu}  =\sum_{i=1}^3 \frac{m_i -1}{2} = p-1.
\end{align}
On the other hand, we shall consider the case (b).
Let us assume, without loss of generality, that $m_1$ is odd and both $m_2$ and $m_2$ are even.
But, since
\begin{align}
0  \leq  -m_1^{\mu} +  m_2^{\mu} +  m_3^{\mu}     = \frac{2p -1 - \sum_{i=1}^3m_i}{2}  = -1, 
\end{align}
this is a constradiction.
Hence, 
 the equality $m =1$ holds.
 One verifies that $m_i \neq 0$ for any $i$.
 Indeed,  one of them, say $m_1$, coincides with $0$ and $m_2$ is odd (resp., even), then 
  \begin{align}
0 \leq -m_1^{\mu} +  m_2^{\mu} +  m_3^{\mu} = - \frac{p-1}{2} + \frac{p-m_2-1}{2} + \frac{p-m_3 -1}{2}  = -1 \\
 \left(\text{resp.,} \  p-2 \geq \sum_{i=1}^3 m_i^{\mu} = \frac{p-1}{2} + \frac{m_2-1}{2} + \frac{m_3-1}{2} = p-1 \right). \notag
 \end{align}
 This is a contradiction.
  Consequently, 
$(m_1, m_2, m_3)$ specifies a strict $p$-branch numbering on $\Gamma^+_\mfP$.
This completes the proof of the lemma.
\end{proof}

Denote by $\overline{\mfC} \mfo_\mfP^{\psi =0}$ the set of (logarithmic) $k$-connections on $\Omega_{\mbP^{1 \mr{log}}/k}$ with vanishing $p$-curvature. (Hence, $\overline{\mfT} \mfa \mfn_\mfP$ specifies a subset of $\overline{\mfC} \mfo_\mfP^{\psi =0}$.) 
Let $(m_1, m_2, m_3)$ be a triple in $p\text{-}\mfB \mfr_{\Gamma^+_\mfp}$ and let $m$ be the integer with $m p = \sum_{i=1}^3 m_i -1$.
Denote by $\mcO_{\mbP^{1 (1)}} (-m)$ a unique (up to isomorphism) line bundle of degree $-m$ on the Frobenius twist $\mbP^{1 (1)}$ of $\mbP^1$.
We obtain the pull-back $F^*(\mcO_{\mbP^{1 (1)}} (-m))$ of $\mcO_{\mbP^{1 (1)}} (-m)$ via the relative Frobenius morphism  $F : \mbP^1 \migi \mbP^{1 (1)}$.
There exists a  $k$-connection $\nabla^\mr{can}$ on $F^*(\mcO_{\mbP^{1 (1)}} (-m))$ with vanishing $p$-curvature  determined uniquely by the condition that the sections of the subsheaf $F^{-1} (\mcO_{\mbP^{1 (1)}} (-m))$ ($\subseteq F^* (\mcO_{\mbP^{1 (1)}} (-m))$) are contained in $\mr{Ker} (\nabla^\mr{can})$ (cf. ~\cite{Wak7}, \S\,1.7).
Also, one may construct uniquely a  $k$-connection $\nabla^{\mr{can}}_{m_1, m_2, m_3}$ on $F^* (\mcO_{\mbP^{1 (1)}} (-m)) (\sum_{i=1}^3 m_i \sigma_i)$ whose restriction to $F^* (\mcO_{\mbP^{1 (1)}} (-m))$ coincides with $\nabla^{\mr{can}}$.
The monodromy of $\nabla^{\mr{can}}_{m_1, m_2, m_3}$ at $\sigma_i$ ($i =1,2,3$) is $-\tau (m_i)$.
Since 
\begin{align}
\mr{deg} (F^* (\mcO_{\mbP^{1 (1)}} (-m)) (\sum_{i=1}^3 m_i \sigma_i)) = -p m + \sum_{i=1}^3 m_i=1 = \mr{deg} (\Omega_{\mbP^{1(1)}/k}),
\end{align}
we have an isomorphism $F^* (\mcO_{\mbP^{1 (1)}} (-m)) (\sum_{i=1}^3 m_i \sigma_i) \isom \Omega_{\mbP^{1(1)}/k}$.
$\nabla^{\mr{can}}_{m_1, m_2, m_3}$ corresponds, via this isomorphism, to a $k$-connection $\breve{\nabla}^{\mr{can}}_{m_1, m_2, m_3}$  on  $\Omega_{\mbP^{1(1)}/k}$ of monodromy $(-\tau (m_1), -\tau (m_2), -\tau (m_3)) \in k^{\times 3}$ (with vanishing $p$-curvature).
Notice that this connection does not depend on the choice of the isomorphism $F^* (\mcO_{\mbP^{1 (1)}} (-m)) (\sum_{i=1}^3 m_i \sigma_i) \isom \Omega_{\mbP^{1(1)}/k}$.
The resulting assignment $(m_1, m_2, m_3)$ $\mapsto \breve{\nabla}^{\mr{can}}_{m_1, m_2, m_3}$ determines a bijection 
\begin{align} \label{E034}
p\text{-}\mfB \mfr_{\Gamma^+_\mfP} \isom \overline{\mfC} \mfo_\mfP^{\psi =0}.
\end{align}
Indeed, its inverse is given by $\nabla \mapsto  (\tau^{-1}(-\mu_1^\nabla), \tau^{-1}(-\mu_2^\nabla), \tau^{-1}(-\mu_3^\nabla))$.
In particular, a $k$-connection on $\Omega_{\mbP^{1 \mr{log}}/k}$ with vanishing $p$-curvature may be uniquely determined by its monodromy.

\vspace{3mm}
\ble \label{P04} \leavevmode\\
 \ \ \ 
The bijection (\ref{E034}) restricts to a bijection
\begin{align} \label{EEr8}
 p\text{-} \mfB \mfr_{\Gamma_\mfP^+}^\mr{st} \isom \overline{\mfT} \mfa \mfn_{\mfP}  
\end{align}
which   makes the following diagram commute:
\begin{align} \label{RR2}
\xymatrix{
 p\text{-} \mfB \mfr_{\Gamma_\mfP^+}^\mr{st} 
  \ar[r]^{(\ref{EEr8})}_{\sim} \ar[d]_{\mu^\mr{comb}_{\Gamma^+}}
&   \overline{\mfT} \mfa \mfn_{\mfP}   
   \ar[r]_{\sim}^{(\ref{EErt6})} & \mfM \overline{\mfO} \mfp^{^\mr{Zzz...}}_{\mfs \mfl_2, \mfP}  \ar[d]^{\mu^{^\mr{Zzz...}}_{\mfP}}  \\
p\text{-} \mfE \mfd_{\Gamma_\mfP^+}^\mr{bal}  \ar[rr]^{\sim}_{(\ref{eq2233})} && \overline{\mfO} \mfp^{^\mr{Zzz...}}_{\mfs \mfl_2, \mfP}.  
}
\end{align}
In particular,
the set $\mfM \overline{\mfO} \mfp^{^\mr{Zzz...}}_{\mfs \mfl_2, \mfP}$ is in bijection with 
the set of triples $(m_1, m_2, m_3)$ consisting of positive integers with $\sum_{i=1}^3 m_i = p+1$.
\ele
\begin{proof}
It is immediately verified that the following diagram is commutative:
\begin{align} \label{RR44}
\xymatrix{
 p\text{-} \mfB \mfr_{\Gamma_\mfP^+}
 \ar[r] \ar[d]
&  \overline{\mfC} \mfo_\mfP
   \ar[r]_{\sim}^{(\ref{GGGG})} & \mfM \overline{\mfO} \mfp_{\mfs \mfl_2, \mfP}  \ar[d]^{\mu_{\mfP}}  \\
\widetilde{\mbF}_p^{\times 3}  \ar[rr]&& \overline{\mfO} \mfp_{\mfs \mfl_2, \mfP},  
}
\end{align}
where
\begin{itemize}
\item[$\bullet$]
 the left-hand vertical arrow denotes the map given by $(m_i)_{i=1}^3 \mapsto ( m_i^{\mu})_{i=1}^3$;
\item[$\bullet$]
 the upper left-hand horizontal arrow is the composite of (\ref{E034}) and the natural injection
$ \overline{\mfC} \mfo_\mfP^{\psi =0} \migiincl  \overline{\mfC} \mfo_\mfP$;
\item[$\bullet$]
the lower horizontal arrow denotes the composite $\widetilde{\mbF}_p^{\times 3} \xrightarrow{\tau^{\times 3}}\mbF_p^{\times 3} \migiincl k^{\times 3} \xrightarrow{(\ref{aaa6})}  \overline{\mfO} \mfp_{\mfs \mfl_2, \mfP}$.
\end{itemize}
The diagram (\ref{RR2}) may be obtained from (\ref{RR44}) by restricting $\widetilde{\mbF}_p^{\times 3}$ and $ \overline{\mfO} \mfp_{\mfs \mfl_2, \mfP}$ to $p\text{-} \mfE \mfd_{\Gamma_\mfP^+}^\mr{bal}$ and $\overline{\mfO} \mfp^{^\mr{Zzz...}}_{\mfs \mfl_2, \mfP}$ respectively.
This completes the proof of the lemma.
\end{proof}
\vspace{3mm}

By applying Proposition  \ref{yy0176}, one may glue together the isomorphisms (\ref{EEr8}) applied to $\mfP_v$ for the various  $v \in V_\mfX$ (with the notation following Lemma \ref{LLL001}).
The resulting  isomorphism,  as displayed in (\ref{EErh8d}) below,
  gives  a combinatorial description of pre-Tango structures (as well as dormant generic Miura $\mfs \mfl_2$-opers)  on an arbitrary   totally degenerate curve.
(Theorem A follows from the following assertion together with the decompositions (\ref{EE002}), (\ref{TTTT}), (\ref{YYYY}), \ref{UUUU}.)

\vspace{3mm}
\bco \label{PP005} \leavevmode\\
 \ \ \ Let $\mfX$ be a  totally degenerate pointed stable curve over $k$ of type $(g,r)$ and  $\vec{\varepsilon}$  an element of $\mbF_p^{\times r}$, where $\vec{\varepsilon} :=  \emptyset$ if $r =0$.
 Then, 
 there exists a canonical bijection
 \begin{align} \label{EErh8d}
\overline{\mfT} \mfa \mfn_{\mfX, \vec{\varepsilon}} \isom  p\text{-} \mfB \mfr_{\Gamma_\mfX^+, -\vec{\varepsilon}}^\mr{st}
\end{align}
which   makes the following diagram commute:
\begin{align} \label{RR003}
\xymatrix{
 \mfM \overline{\mfO} \mfp^{^\mr{Zzz...}}_{\mfs \mfl_2, \mfX, \vec{\varepsilon}\cdot \check{\rho}}   \ar[d]_{\mu^{^\mr{Zzz...}}_{\mfs \mfl_2, \mfP}} & \overline{\mfT} \mfa \mfn_{\mfX, \vec{\varepsilon}} \ar[r]_{\sim}^{(\ref{EErh8d})}  \ar[l]^{\sim}_{(\ref{EE007})}  & p\text{-} \mfB \mfr_{\Gamma_\mfX^+, -\vec{\varepsilon}}^\mr{st} \ar[d]^{\mu^\mr{comb}_{\Gamma_\mfX^+, -\vec{\varepsilon}}} \\
\overline{\mfO} \mfp^{^\mr{Zzz...}}_{\mfs \mfl_2, \mfX, \chi (\vec{\varepsilon}\cdot \check{\rho})}   \ar[rr]^{\sim}_{(\ref{eq8901})} &  & p\text{-} \mfE \mfd_{\Gamma_\mfX^+, -\vec{\varepsilon}^{\, \mu}}^\mr{bal}.
}
\end{align}
In particular, the set $\mfM \overline{\mfO} \mfp^{^\mr{Zzz...}}_{\mfs \mfl_2, \mfX,  \vec{\varepsilon}\cdot \check{\rho}}$ is in bijection with the set $p\text{-} \mfB \mfr_{\Gamma_\mfX^+, -\vec{\varepsilon}}^\mr{st}$.
\eco

\vspace{5mm}
\subsection{Case of $g =1$} \label{SSErt}
\leavevmode\\ \vspace{-4mm}

Finally, we shall conclude the present paper with proving 
 that there is no strict $p$-branch numbering on $\Gamma^+$
 unless the underlying semi-graph $\Gamma$ is of type $(g,r)$ with $g \leq 1$ (cf. Proposition \ref{P048} below).
Equivalently, there is no dormant generic Miura $\mfs \mfl_2$-oper on any totally degenerate  curve  of genus $g >1$ (cf. Corollary \ref{PCD04}).

\bpr \label{P048} \leavevmode\\
 \ \ \ 
 Let $\Gamma^+ : = (V_\Gamma, E_\Gamma, \zeta_\Gamma, \lambda_\Gamma)$ be a marked semi-graph whose underlying semi-graph $\Gamma$ is connected,  $3$-regular, and of type $(g,r)$.
\begin{itemize}
\item[(i)]
If there exists a strict $p$-branch numbering on $\Gamma^+$, then
   the inequality $g \leq 1$ holds.
\item[(ii)]
 If $g =1$, then the  inclusion
\begin{align}
  p\text{-} \mfB \mfr_{\Gamma^+, \vec{e}}^\mr{st} \migiincl p\text{-} \mfB \mfr_{\Gamma^+}^\mr{st}
\end{align} 
is bijective, where $\vec{e} := (\overline{-1}, \overline{-1}, \cdots, \overline{-1}) \in \mbF_p^{\times r}$ if $r >0$ (resp., $\vec{e} := \emptyset$ if $r =0$),
and the following equalities hold:
\begin{align}
\sharp (  p\text{-} \mfB \mfr_{\vec{e}, \Gamma^+}^\mr{st}) = \sharp (p\text{-} \mfB \mfr_{\Gamma^+}^\mr{st}) =p-1.
\end{align}
 \end{itemize}
 \epr
\begin{proof}
First, we shall consider assertion (i).
Let $\vec{m} := (m_b)_{b \in  B_\Gamma} \in \widetilde{\mbF}_p^{B_\Gamma}$ be  a strict $p$-branch numbering  on $\Gamma^+$. 
Let us assume that $g >0$.
Since $H_1 (\Gamma, \mbZ) \neq 0$ (where we  regard $\Gamma$ as a topological space in the manner  mentioned  in Remark \ref{R0ff01}), one may find a vertex $v_0$ of $\Gamma$ and a path $(b_{j})_{j =1}^l$  from $v_0$ to $v_0$ itself 
such that $b_j \neq b_{j-1}^\divideontimes$ for any $j \in \{ 1, \cdots, l \}$, where $b_0^\divideontimes := b_l^\divideontimes$.
(We shall refer to such a path as a {\bf reduced loop} based at $v_0$.)
For each $j \in \{1, \cdots, l \}$, there exists a unique branch $b^{c}_j \in B_\Gamma$ with  $B_{\zeta_\Gamma (b_j)} := \{ b_j, b_{j-1}^\divideontimes, b_j^c \}$.
The assumption  that $\vec{m}$ specifies a strict $p$-branch numbering implies the following equalities
\begin{align}
m_{b_1} & = p+1 -(m_{b_0^\divideontimes} + m_{b_1^c}), \\
m_{b_2} & = p+1 -(m_{b_1^\divideontimes} + m_{b_2^c}) \notag \\
& = p+1 -((p- m_{b_1}) + m_{b_2^c}) \notag \\
& = p+2 -(m_{b_0^\divideontimes} + m_{b_1^c} + m_{b_2^c}),\notag \\
m_{b_3} & = \cdots \notag \\
\vdots & \ \  \vdots \  \  \vdots \notag \\
m_{b_l} & = p + l -(m_{b_0^\divideontimes} + \sum_{j=1}^l m_{b_j^c})  \notag \\
& = l  + m_{b_l} - \sum_{j=1}^l  m_{b_j^c}. \notag
\end{align}
The last equality is equivalent to the equality $l = \sum_{j=1}^l  m_{b_j^c}$, which implies the equality $m_{b_j^c} =1$
for any $j \in \{1, \cdots, l \}$.

Now,  suppose further that $g \geq 2$.
Then, one verifies from  the topological structure of $\Gamma$  that (after replacing $v_0$ by another vertex)
there exist
two reduced loops $(b'_{j})_{j=1}^{l'}$, $(b''_{j})_{j=1}^{l''}$ based at $v_0$ such that
$\{ b'_1, {b'}_l^\divideontimes \} = \{ b_l^\divideontimes, b_1^c \}$, $\{ b''_1, {b''}_l^\divideontimes \} = \{ b_1^c, b_1 \}$.
By applying the above discussion to $(b'_{j})_{j=1}^{l'}$ and  $(b''_{j})_{j=1}^{l''}$ respectively,  
we obtain the equalities $(m_{b_1^c} =) \ m_{b_1} = m_{b_l^\divideontimes} = 1$.
Hence,  $3 \ (= m_{b_1} + m_{b_1^c} + m_{b_l^\divideontimes}) = 1 +p$,   and this contradicts the assumption  $p >2$.
Consequently, the inequality $g \leq 1$ holds.

Next, let us consider 
 assertion (ii).
Suppose that  $g =1$, and let $v_0$  and $(b_j)_{j=1}^l$ be as above.
If we are given a strict $p$-branch numbering $\vec{m} := (m_b)_{b \in B_\Gamma}$ on $\Gamma^+$,
then  since $m_{b_j^c} =1$ and $m_{b_j} + m_{b_j^c} + m_{b_{j-1}^\divideontimes} = 1 +p$ (for any $j$),
there exists a unique $a \in \{1, \cdots, p-1 \}$ satisfying the condition $(*)_a$ described as follows: $(*)_a$   $m_{b_j} = a$ and $m_{b_j^\divideontimes} = p-a$ for any $j$.
Conversely, for each $a \in \{ 1, \cdots, p-1 \}$, one may construct a unique   strict $p$-branch numbering $\vec{m} := (m_b)_{b \in B_\Gamma}$ on $\Gamma^+$ satisfying the condition 
$(*)_a$
 in such a way that for each $v \in B_\Gamma \setminus \{ b_j\}_{j=1}^l$, the multiset $[m_{b_1}, m_{b_2}, m_{b_3}]$ (where $B_v = \{b_1, b_2, b_3 \}$) coincides with $[1, 1, p-1]$. 
If $r >0$, then each such   strict $p$-branch numbering is verified to be of exponent $\vec{e} \in \mbF_p^{\times r}$.
Thus, this completes the proof of  assertion (ii).
\end{proof}
\vspace{3mm}

\bco \label{PCD04} \leavevmode\\
 \ \ \ 
 Let $\mfX$ be a  totally degenerate curve over $k$ of type $(g,r)$.
 If there exists a dormant generic Miura $\mfs \mfl_2$-oper $\widehat{\mcE}^\spadesuit$ on $\mfX$, 
 then the inequality $g \leq 1$ holds.
 If, moreover,  $g =1$ (and $r >0$), then there exist  precisely $p-1$ dormant generic Mura $\mfs \mfl_2$-opers on $\mfX$,  and  these are 
  of exponent $\vec{e}$ (cf. Proposition \ref{P048} (ii) for the definition of $\vec{e}$).
 \eco
\begin{proof}
The assertion follows from Corollary \ref{PP005}  and Proposition \ref{P048}.
\end{proof}
\vspace{3mm}

\end{document}